\documentclass[11pt,a4paper,leqno]{article}
\usepackage{amsfonts,amssymb,mathrsfs,amsmath}

\usepackage{amscd}

\vspace{0.5cm}

 4

\font\erm=cmr8

\usepackage{color}
\usepackage{graphicx}

\author{Andrzej Krzysztof Kwa\'sniewski}
\title{Note on Ward-Horadam \textit{H(x)}- binomials' recurrences and related interpretations, II}

\newtheorem{defn}{Definition}

\newtheorem{theoremn}{Theorem}

\newtheorem{quot}{Quotation}

\newcommand{\fnomial}[2]{ {{#1} \choose {#2}}_{\!\!F} }
\newcommand{\fnomialF}[3]{ {{#1} \choose {#2}}_{\!\!#3} }

\newcommand{\layer}[2]{\langle\Phi_#1 \to \Phi_#2\rangle}

\begin{document}

\begin{center}
\noindent {\sffamily Note on Ward-Horadam \textit{H(x)}- binomials' recurrences and related interpretations, II}  \\ 

\vspace{0.5cm} 
\sffamily Andrzej Krzysztof Kwa\'sniewski

\vspace{0.5cm}
{\erm
\sffamily Member of the Institute of Combinatorics and its Applications, 
Winnipeg, Manitoba, Canada \\
	PL-15-674 Bia\l ystok, Konwaliowa 11/11, Poland\\
	e-mail: kwandr@gmail.com\\}
\end{center}

\vspace{0,5cm}

\noindent \textbf{Summary} 

\noindent \sffamily This note is a continuation of \cite[2010]{akk2010}. \textbf{Firstly}, we propose \textit{H(x)}-binomials' recurrence formula appointed by Ward-Horadam $ H(x)= \left\langle H_n(x)\right\rangle _{n\geq 0}$ functions' sequence i.e. any functions' sequence solution of the second order recurrence with functions' coefficients. As a method this  comprises    $H \equiv H(x=1)$ number sequences \textit{V}-binomials' recurrence formula determined by the primordial  Lucas sequence of the second kind   $V = \left\langle V_n\right\rangle_{n\geq 0}$ as well  as its well elaborated companion fundamental Lucas sequence of the first kind $ U = \left\langle U_n \right\rangle_{n\geq 0}$ which gives rise in its turn to the known  \textit{U}-binomials' recurrence as in \cite[1878]{EdL} , \cite[1949]{JM},  \cite[1964]{TF},  \cite[1969]{Gould},  \cite[1989]{G-V}  or in  \cite[1989]{K-W}  etc. 

\noindent Then we deliver a new type \textit{H(x)}-binomials' "mixed" recurrence,  straightforward formula for special case of Ward-Horadam sequences $H(x) = \left\langle p^n(x) + q^n(x)\right\rangle_{n\geq 0}$.

\vspace{0,1cm}

\noindent Here general $H(x)$-binomials' array is  appointed by Ward-Horadam sequence of functions which in predominantly considered cases where chosen to be polynomials.

\vspace{0,1cm} 

\noindent For the sake of combinatorial interpretations and in number theory  $H(x=1) = \left\langle H_n(x=1) \right\rangle_{n\geq 0}$ is  usually considered to be natural or integer numbers valued sequence. Number sequences $H =H(x=1)= \left\langle H_n\right\rangle _{n\geq 0}$ were also recently called by several authors: Horadam sequences.

\vspace{0,2cm}

\noindent \textbf{Secondly}, we supply a review  of  selected related combinatorial interpretations of generalized binomial coefficients.  We then propose also a kind of transfer of  interpretation of $p,q$-binomial  coefficients onto $q$-binomial interpretations thus bringing us back to Gy{\"{o}}rgy P\'olya  and Donald Ervin  Knuth relevant investigation decades ago.

\vspace{0,1cm}

\noindent The list of references is prevailingly  indicatory (see references therein) and is far from being complete.

\vspace{0,2cm}
\noindent AMS Classification Numbers: 05A10 , 05A30.

\vspace{0,1cm}

\vspace{0.4cm}


\section{General Introduction}

\noindent \textbf{1.1.  $p,q$ people are followers of Lucas people.} The are many authors who use in their investigation the fundamental Lucas sequence $ U \equiv \left\langle n_{p,q}\right\rangle_{n\geq 0}$ - frequently  with different notations - where $n_{p,q} = \sum_{j=0}^{n-1}{p^{n-j-1}q^j}=U_n$; see Definition 1  and then definitions that follow it.  In regard to this a brief intimation   is on the way.

\vspace{0,1cm}

\noindent Up to our knowledge it was Fran\c{c}ois \'Edouard Anatole Lucas in \cite[1878]{EdL} who was the first who had  \textit{not only} defined \textit{fibonomial} coefficients as stated in \cite[1989]{K-W}   by Donald  Ervin  Knuth and Herbert Saul Wilf but who was the first who had defined  \textit{$U_n \equiv n_{p,q}$}-binomial coefficients ${n \choose k}_U \equiv {n \choose k}_{p,q}$  and had derived a recurrence for them: see page 27, formula (58) \cite[1878]{EdL}. Then - referring to Lucas - the investigation relative to  divisibility properties of relevant number Lucas sequences $D$, $S$ as well as numbers' $D$ - binomials and  numbers' $D$ - multinomials was continued in \cite[1913]{Carmichel} by  Robert Daniel Carmichel;  see pp. 30,35 and 40 in \cite[1913]{Carmichel}  for $ U \equiv D = \left\langle D_n\right\rangle_{n\geq 0}$ and ${n \choose {k_1,k_2,...,k_s}}_D$ - respectively. Note there also formulas (10), (11) and (13) which might perhaps  serve to derive explicit untangled form of recurrence for the  $ \textit{V}$- \textit{binomial coefficients} ${n \choose k}_V \equiv {n \choose k}_S$ denoted by  primordial Lucas sequence $\left\langle S_n\right\rangle_{n\geq 0}= S \equiv V$. Let $F(x=1)\equiv A$ denotes a number sequence  as in \cite[1915]{Fon}  by Fonten\'e Georges, \cite[1969]{Gould}  by Henri W.Gould  then followed by for example
\cite[2005]{ST 2005},  \cite[2005]{ST 2005 2}, \cite[2006]{ST 2006} by Jaroslav Seibert  and Pavel Trojovsk\'y  and  \cite[2007]{T 2007}  by  Pavel Trojovsk\'y.

\vspace{0,1cm}

\noindent $A$ - multinomial coefficients' \textit{recurrences} are not present in that early and other works and up  to our knowledge a special case of such appeared at first in \cite[1979]{Shannon 1979 multi} by  Anthony G. Shannon. More on that - in what follows after Definition 3.

\vspace{0,1cm}

\noindent Significant peculiarity of Lucas originated sequences  includes  their importance for number theory  (see  middle-century paper \cite{Halton 1996} by  John H. Halton and recent, this century papers \cite[2010]{Smyth1} by Chris Smith  and \cite[2010]{Smyth2} by K\'alm\'an  Gy{\"{o}}ry with Chris Smith  and the reader may enjoy also the PhD Thesis  \cite[1999]{Anne-Marie 1999}  by  Anne-Marie Decaillot-Laulagnet). This Lucas originated investigation amalgamates diverse areas of mathematics due to hyberbolic - trigonometric character of these Fonctions Num\'eriques Simplement Périodiques i.e. fundamental and primordial Lucas sequences - as beheld  in \cite[1878]{EdL}. One may track then  a piece of further constructions  for example in \cite[1999]{akk1999}). \\

\noindent There in \cite[1999]{akk1999} tail formulas (3.12) and (3.14) are illustrating   the proved and exploited by \'Eduard Lucas  complete analogy of the $V_n$  and $U_n$  symmetric functions of roots with the circular and hyperbolic functions of order $ 2$. This is so  due to Lucas formulas (5) in \cite{EdL} rewritten in terms of $\cosh$  and $\sinh$ functions as formulas (3.13) and (3.14) in \cite{akk1999} where there these  result from   de Moivre one parameter group  introduced in \cite{akk1999}  in order to pack compactly  the addition formulas (1.6), (1.7) in \cite{akk1999} into abelian group  "parcel"   encompassing  Pafnuty Lvovich Chebyshev polynomials of both kinds. 

\vspace{0,1cm}

\noindent In this connection see the Section 2 in the  recent Ward-Horadam people paper  \cite [2009]{He-Siue}  by  Tian-Xiao He,  Peter Jau-Shyong Shiue. There in Proposition 2.7.  illustrative   Example 2.8. with Pafnuty Lvovich Chebyshev polynomials of the first kind the well known recurrence formula (2.28)  is equivalent to  abelian one-parameter de Moivre matrix group  multiplication rule from which the corresponding recurrence  (1.7) in  \cite[1999]{akk1999}  follows.

\vspace{0,1cm}


\noindent \textbf{1.2. We deliver here}  - continuing  the note \cite{akk2010} - a new  \textit{H(x)}-binomials' recurrence formula appointed by Ward-Horadam $ H(x)= \left\langle H_n(x)\right\rangle _{n\geq 0}$ field of zero characteristic nonzero valued functions' sequence which comprises  for $H \equiv H(x=1)$ number sequences case - the  \textit{V}-binomials' recurrence formula determined by the primordial Lucas sequence of the second kind  $V = \left\langle V_n\right\rangle_{n\geq 0}$   \cite[2010]{akk2010} as well  as its well elaborated companion fundamental Lucas sequence of the first kind $ U = \left\langle U_n\right\rangle_{n\geq 0}$ which gives rise in its turn to the  \textit{U}-binomials' recurrence as in \cite[1878]{EdL} , \cite[1949]{JM},  \cite[1964]{TF},  \cite[1969]{Gould},  \cite[1989]{G-V}  or in  \cite[1989]{K-W}  and so on. 

\vspace{0,1cm}

\noindent We do it by following recent applicable work \cite[2009]{Savage}  by Nicolas A. Loehr and Carla  D.  Savage   thought one may  - for that purpose - envisage now easy extensions of particular $p,q$ - cases considered earlier - as  \textit{for example} the following: the relevant recursions in \cite[1949]{JM}, in  \cite[1989]{G-V},  in \cite[1992]{KaKi} - ( recursions  (40) and (51)) , or \cite[2000]{Holte 2000} by  John  M.  Holte (Lemmas 1,2 dealing with  $U$-binomials provide a  motivated example for observation  Theorem 17 in \cite{Savage} ) One is invited also to track Lemma 1 in  \cite[2001]{Hu-Sun} by  Hong Hu  and   Zhi-Wei Sun ; see also corresponding recurrences for $p,q$-binomials $\equiv$ $U$-binomials in \cite[1878]{EdL} or in \cite[2008]{MD2} \textbf{v[1]} by Maciej Dziemia\'nczuk (compare there (1) and (2) formulas), or see Theorem 1  in \cite[2008]{Corsino}  by Roberto Bagsarsa Corcino as well as  track the proof of the Corollary 3. in \cite[2009]{MD2} \textbf{v[2]}  by Maciej Dziemia\'nczuk.

\vspace{0,1cm}


\noindent \textbf{This looked for} here new \textit{H(x)}-binomials' overall recurrence formulas (recall: encompassing \textit{V}-binomials for primordial Lucas sequence $V$)  is not present neither in \cite{EdL} nor in \cite{Savage}, nor in \cite[1915]{Fon} , nor in \cite[1936]{Ward}, nor in \cite[1949]{JM}. Neither we find it in  - quoted here contractually by a nickname as  "Lucas $(p,q)$-people" - references [1-47]. Neither it is present - up to this note author knowledge - in  all other   quoted here contractually by a nickname as  "Ward-Horadam -people" - references [53-84]. As for  "Lucas $(p,q)$-people" and "Ward-Horadam -people" references - (including these [n]  with  $n>86$  - the distinction  which People are which is  quite contractual. The nicknames are nevertheless indicatively helpful.  We shall be more precise soon - right with definitions are being started.\\
\noindent \textbf{Interrogation}. Might it be so that such an  overall recurrence formula does not exist? The method to obtain formulas is more general than actually  produced family of formulas?  See Theorems 2a and 2b, and Examples 4,2, ... in what follows.


\vspace{0,1cm}

\noindent \textbf{Meanwhile}  \textit{H(x)}-binomials' overall recurrence formula for the  Ward-Horadam's function sequence  $H(x) = \left\langle H_n(x)\right\rangle _{n\geq 0}$  follows straightforwardly from the easily proved important observation - the Theorem 17 in \cite[2009]{Savage} as already had it been remarked in \cite[2010]{akk2010} for the $H \equiv H(x=1)$ case. Finally, see this note observation named Theorem 2b.

\vspace{0,1cm}

\noindent \textbf{This paper formulas} may and should be confronted with  Fonten\'e obvious recurrence for complex valued $A$-binomials ${n \choose k}_A$,  $A \equiv A(x=1)$  in \cite[1915]{Fon} i.e. with (6) or (7) identities in \cite[1969]{Gould} by Henri W. Gould or with recurrence in \cite[1999]{ALs} by Alexandru Ioan  Lupas ,  which particularly also stem easily just  from the definition of \textbf{any} $F(x)$-binomial coefficients arrays with $F(x) = \left\langle F_n(x)\right\rangle_{n\geq 0}$ staying for any field of characteristic zero  nonzero valued functions'sequence ; $F_n(x) \neq 0,\:n\geq 0$.  For $F = F(x=1)$-multinomial coefficients automatic definition see \cite[1913]{Carmichel}  by    Robert Daniel Carmichel or then  \cite[1969]{Gould} by Henri W. Gold and finally  see \cite[1979]{Shannon 1979 multi} by Anthony G. Shannon, where recurrence is proved for ${n \choose {k_1,k_2,...,k_s}}_U$    with $U$-Lucas fundamental being here complex valued number sequence. For $F(x)$ - multinomial coefficients see  \cite[2004]{akk3 2004}  and compare with  $F(x)$-binomials from \cite[1999]{ALs} or those from \cite[2001]{Rich 2001}.

\vspace{0,1cm}

\noindent To this end we supply now two informations (1.3)  and  (1.4) pertinent  ad references and ad nomenclature.

\vspace{0,2cm}

\noindent \textbf{1.3. Ad the number theory and divisibility properties references.}  
\noindent For the sake of combinatorial interpretations of $F$- number sequences as well as their correspondent $F$-multinomial coefficients  and also for the sake of the number theoretic studies of Charles Hermite  \cite{Herm 1889-1} and with Thomas Jan Stieltjes in   \cite{Herm 1889-2} or by  Robert Daniel Carmichel  \cite[1913]{Carmichel} or  \cite[1919]{Dick} or that of Morgan Ward  \cite[1936]{Ward1}, \cite[1939]{Ward2}, \cite[1937]{Ward3}, \cite[1937]{Ward4}, \cite[1954]{Ward0}, \cite[1955]{Ward5}, \cite[1959]{Ward6} and that of Derric Henry Lehmer  \cite[1930]{Leh}, \cite[1933]{Leh1},  \cite[1935]{Leh2} or this of  Andrzej Bobola Maria Schinzel   \cite[1974]{Schinzel} and Others' studies  \textit{on  divisibility properties} -  these are the sub-cases  $F_n \in\mathbb{N}$  or $F_n \in\mathbb{Z}$ which are being regularly considered at the purpose.  

\vspace{0,1cm}

\noindent   As for the "Others" - see \textit{for example}:  \cite[1959]{Durst}, \cite[1973]{RAKKK},  \cite[1974]{SingMa1},  \cite[1974]{SingMa2},  \cite[1974]{SingMa3},  \cite[1973]{SingMa4}, \cite[1977]{KKK},  \cite[1977]{Stewart}, \cite[1979]{ShaHor},  \cite[1979]{HorLohSha 1979},  \cite[1980]{Somer},  \cite[1980]{Beukers},  \cite[1989]{K-W}, \cite[1991]{SatoShi}, \cite[1992]{Zhi-Wei Sun}, \cite[1995]{GouldPaula},  \cite[1999]{adGould},  \cite{GouldList},   \cite[1995]{Hil-Vranc}, \cite[1995]{KimWeb}, \cite[1995]{Voutier},  \cite[1998]{Wilson},    \cite[2001]{Hu-Sun}, \cite[2006]{MoAbo}, \cite[2009]{MDWB}. 

\vspace{0.2cm}

\noindent \textbf{1.4. Ad the name: \textcolor{blue}{Ward-Horadam sequence}.} 
\noindent According to the authors of  \cite[2009]{He-Siue}  it was  Mansour   \cite{Mansour}  who called  the  sequence $H = \left\langle a_n \right\rangle_{n\geq 0}$  defined by (1)  a \textbf{Horadam's sequence}, as - accordingly to the author of \cite{Mansour} - the number sequence $H$ was introduced in 1965 by Horadam  \cite{Horadam1} (for special case of Ward-Horadam number sequences see Section 2 in \cite[1974]{WaltHorad} and see also  \cite[2009]{HorzumKocer}), this however  notwithstanding  the ingress of complex numbers valued $F$-\textit{binomials}  and $F$-\textit{multinomials} into Morgan Ward's systematic \textit{Calculus of sequences} in  \cite[1936]{Ward} and then in 1954 Ward's introduction   of "`\textit{nomen omen}"'  $W\equiv H$ in \cite[1954]{Ward0} integer valued sequences. 

\vspace{0,1cm}

\noindent Perceive then the appraisal of adequate Morgan Wards' work in the domain by  Henri W.Gould \cite[1959]{Gould} and by Alwyn F. Horadam   and Anthony G. Shannon  in  \cite[1976]{HorSha}  or Derrick  Henry Lehmer in \cite[1993]{Leh3}. 
\noindent On this occasion note also the \textbf{Ward-Horadam} number sequences \textbf{in} \cite[1965]{Zeitlin0}  and  \cite[1965]{Zeitlin}.    

\vspace{0,1cm}

\noindent The  sequence $H = \left\langle a_n \right\rangle_{n\geq 0}$  defined by (2)  was called   \textbf{Horadam's sequence}  in \cite[1974]{Hilton 1974 partition}  by Anthony J. W. Hilton and in generalization of \cite[1974]{Hilton 1974 partition} which is the work \cite[1979]{Shannon 1979 Hilton's partition}  by  Anthony G. Shannon where consequences of the partition of the set of nontrivial  solutions of (1) into solutions of Lucas type $L$ (primordial Lucas sequences)   and solutions of Fibonacci type $F$ (fundamental Lucas sequences) were studied by both authors. To be more precise: the article \cite[1974]{Hilton 1974 partition} opens the way to consider $F$-binomial and $L$-binomial coefficients (not considered in \cite[1974]{Hilton 1974 partition})  while the presentation  \cite[1979]{Shannon 1979 Hilton's partition} opens the way to consider $V^{(r)}_s$-$r$-nomial coefficients where $s = 1,...,r$ (neither considered in \cite[1979]{Shannon 1979 Hilton's partition}).  It should be though  forthwith noted that  in another paper \cite[1979]{Shannon 1979 multi} Anthony G. Shannon derives  the recursion for  fundamental $V^{(r)}_r $-$r$-nomial coefficients denoted there as $u^{(r)}$-multinomial coefficients, where there specifically  for $r=2$ the number sequence  $u^{(2)}$ coincides with  the fundamental Lucas sequence (i.e. the Lucas sequence of the first kind  $\equiv$  of the $F$-type in the nomenclature of \cite[1974]{Hilton 1974 partition}  by Anthony J. W. Hilton). See more - what follows right after this note Definition 3  and finally compare with the way indicated by this note Theorem 2b.

\vspace{0.2cm}

\noindent \textbf{1.5.} The  sequence $H = \left\langle a_n \right\rangle_{n\geq 0}$  defined by the recurrence (2) - below -  with the initial conditions   $H_0 =1$ and  $H_1 =s$  was exploited in \cite[1962]{Carlitz} by Leonard Carlitz.  Soon,  Leonard Carlitz in  \cite[1965]{Carlitz 1965} had  proved the following $H$-binomial formula for the sequence $H = \left\langle a_n \right\rangle_{n\geq 0}$  defined by the recurrence (2) below  with $t=1$ and \textit{distinct} roots $p,q$ :

$$ \prod_{j=0}^n\left(x-p^j\cdot q^{n-j} \right) = \sum_{r=0}^{n+1}\left( -1 \right)^{\frac{r(r+1}{2}}\cdot {n+1 \choose r}_{U^*}\cdot x^{n+1-r},$$
where
$${n \choose k}_{U^*} \equiv {n \choose k}_U $$
and   $U^*_n = p^n - q^n$  while  $U_n = \frac{p^n - q^n}{p-q}$. Compare with this note Theorem 2a.


\section{Preliminaries}

\vspace{0.1cm}

\textit{Names}:  The Lucas sequence   $V = \left\langle V_n\right\rangle_{n\geq 0}$  is called the Lucas sequence of the second kind - see: \cite[1977, Part I]{KKK}, or \textbf{primordial} - see
\cite[1979]{ShaHor}.\\     

\noindent The Lucas sequence  $ U = \left\langle U_n\right\rangle_{n\geq 0}$   is called the Lucas sequence of the first kind - see: \cite[1977, Part I]{KKK}, or \textbf{fundamental}  - see p. 38 in \cite[1949]{JM} or see \cite[1979]{Shannon 1979 multi}  and  \cite[1979]{ShaHor}.     

\vspace{0.1cm} 

\noindent In the sequel we shall deliver the looked for recurrence for $H$-binomial coefficients ${n \choose k}_H$  determined by the Ward-Horadam sequence $H$ - defined below.

\vspace{0.1cm}

\noindent In compliance with Edouard Lucas' \cite[1878]{EdL} and twenty, twenty first century $p,q$-people's notation we shall at first review here in brief the general second order recurrence;
(compare this review  with the  recent "Ward-Horadam" peoples' paper  \cite [2009]{He-Siue} by Tian-Xiao He and Peter Jau-Shyong Shiue or earlier  $p,q$-papers \cite[2001]{Sun-Hu-Liu} by Zhi-Wei Sun, Hong Hu, J.-X. Liu  and \cite[2001]{Hu-Sun}  by Hong Hu and Zhi-Wei Sun). And with respect to natation:  If in  \cite[1878]{EdL} Fran\c{c}ois  \'Edouard Anatole Lucas  had been used $a=\textbf{p}$ and $b = \textbf{q}$ notation,  he would be perhaps at first glance notified and recognized as a Great Grandfather of all the $(p,q)$ - people. Let us start then introducing reconciling and matched denotations and nomenclature.

\begin{equation}
H_{n+2} = P \cdot H_{n+1}  - Q \cdot H_n ,\ \;n\geq 0 \ and \;  H_0 = a,\; H_1 = b.                                                                                     
\end{equation}
which is sometimes being written  in $\left\langle P,-Q \right\rangle \mapsto  \left\langle s,t \right\rangle$ notation.

\begin{equation}
H_{n+2} = s \cdot H_{n+1}  + t\cdot H_n ,\ \;n\geq 0 \ and \;  H_0 = a,\; H_1 = b.                                                                                     
\end{equation}
We exclude the cases when (2) is recurrence of the first order, therefore we assume that the roots $p,q$ of (5) are distinct   $p\neq q$  and $\frac{p}{q}$ is not the root of unity. We shall come back to this finally while formulating this note  observation named Theorem 2a.

\vspace{0.2cm}
\noindent Simultaneously and collaterally  we  mnemonically pre adjust the starting point to discuss the $F(x)$ polynomials' case via - if entitled - antecedent "$\mapsto$ action": 
$H \mapsto H(x)$, $s \mapsto s(x)$, $t \mapsto t(x)$, etc. 

\begin{equation}
H_{n+2}(x) = s(x) \cdot H_{n+1}(x) + t(x)\cdot H_n, \; n\geq 0, \ H_0 = a(x), H_1 = b(x).                                                                                     
\end{equation}

\vspace{0.2cm}
\noindent enabling recovering  explicit formulas also for sequences of  polynomials correspondingly generated by the above linear recurrence of order 2  - with  Tchebysheff  polynomials and the generalized Gegenbauer-Humbert polynomials included. See for example Proposition 2.7 in the  recent Ward-Horadam peoples' paper  \cite [2009]{He-Siue}  by  Tian-Xiao He and Peter Jau-Shyong Shiue.

\vspace{0.2cm}

\noindent The general solution of (1):   $ H(a,b;P,Q) = \left\langle H_n\right\rangle _{n\geq 0}$ is being called  throughout  this paper - \textbf{Ward-Horadam number'sequence}.

\vspace{0.1cm}

\noindent The general solution $ H(x) \equiv H(a(x),b(x);s(x),t(x)) = \left\langle H_n(x)\right\rangle _{n\geq 0}$ of the recurrence (3) is being called  throughout this paper - \textbf{Ward-Horadam functions' sequence}. It is then to be noted here that ideas germane to \textit{special  Ward-Horadam   polynomials  sequences}  of the \cite{1 Horadam 1996}  paper  were already explored in some details in  \cite{Horadam1}.  For more on special  Ward-Horadam   \textit{polynomials  sequences}  by Alwyn F. Horadam - consult  then:  \cite{Horadam 1969},  \cite[1985] {Horadam3},   \cite{Horadam 1993} ,  \cite{2 Horadam 1996} or  see \textit{for example} the following papers and references therein:   recent papers  \cite[2009]{HorzumKocer} by Tugba Horzum and  Emine  G{\"{o}}kcen Kocer  and  \cite [2009]{Gi-Sang Cheon}  by Gi-Sang Cheon, Hana Kim   and   Louis W. Shapiro.  For   \textit{Ward-Horadam   functions  sequences}    \cite[2009]{He-Siue}  by Tiang-Xiao He and Peter Jau-Shyong Shiue who however there then concentrate on  on \textit{special  Ward-Horadam   polynomials  sequences} only.  Other "polynomial" references  shall appear in the course of further presentation.

\vspace{0.1cm}

\noindent For example - in \cite[2010]{Cigler q-Lucas 2010} Johann Cigler considers special Ward-Horadam   \textit{polynomials  sequences} and  among others he supplies the tiling combinatorial interpretation of these special Ward-Horadam   \textit{polynomials  sequences} which are  $q$-analogues of the Fibonacci and Lucas polynomials introduced in \cite[2002]{Cigler q-Fib q-Luc 2002}  and  \cite[2003]{Cigler Fib Pol 2003} by  Johann Cigler.

\vspace{0.1cm}

\noindent In the paper  \cite[2003]{Cigler} Johann Cigler  introduces "abstract Fibonacci polynomials"  - interpreted in terms of Morse coding sequences monoid with concatenation (monominos and dominos tiling then)   Cigler's abstract Fibonacci polynomial sare monoid algebra over reals valued polynomials with straightforward Morse sequences i.e. tiling recurrence originated (1.6)  "addition formula"

$$F_{m+n}(a,b) =  F_{m+1}(a,b)\cdot F_{m}(a,b) + b\cdot F_{n-1}(a,b)\cdot F_{n}(a,b),$$
which is attractive and seductive to deal with   within the context of this paper Theorem 1 below.                                                                                                                                                                                                                       

\vspace{0.2cm}

\noindent  From the characteristic equation of (1)

\begin{equation}
x^2 = P \cdot x  - Q,                                                                                      
\end{equation}
written by some of $p,q$-people  as

\begin{equation}
x^2 = s \cdot x  + t                                                                                      
\end{equation}
we readily find the Binet form solution of (2) (see (6) in \cite[2009]{HorzumKocer})  which is given by (6) and (7).

\begin{equation}
H_n(a,b;P,Q) \equiv H_n(A,B;p,q) = A p^n  + B q^n, \;n\geq 0 , H_0 = a, H_1 = b.                                                                                     
\end{equation}
where $p,q$ are roots of (5) and we have assumed since now on that  $p\neq q$  and $\frac{p}{q}$ is not the root of unity (see Lemma for example lemma 1 in \cite[1980]{Beukers} 
by Fritz Beukers). Hence for established $p \neq q$  we shall use sometimes the shortcuts
$$H(A,B)\equiv H_n(A,B;p,q) = (p-q)\cdot U_n(A,-B;p,q)\equiv (p-q)\cdot U_n(A,-B). $$

\noindent As for the case $p=q$  included see for example Proposition 2.1 in  \cite[2009]{He-Siue} and see references therein.

\vspace{0.2cm}

\noindent Naturally :  $p+q = P\equiv s$  , $p \cdot q = Q \equiv -t$  and 

\begin{equation}
A = \frac{b-qa}{p-q} \ , \;    B = - \frac{b-pa}{p-q}.   
\end{equation}
\noindent hence we may and we shall use the following  conventional identifications-abbreviations

\begin{eqnarray}
H(A,B) \equiv H(a,b;P,Q) \equiv  H(A,B;P,Q) \equiv  H(A,B;p,q)=	%
\\
=
(p-q)\cdot U_n(A,-B;p,q)\equiv (p-q)\cdot U_n(A,-B).
\nonumber
\end{eqnarray}

\vspace{0.2cm}

\noindent It is obvious that the exponential generating function for Ward-Horadam sequence $H$ reads:

\begin{equation}
E_H(A,B;p,q)[x] = A exp[p\cdot x] + B exp[q\cdot x].   
\end{equation}

\vspace{0.2cm}

\noindent The derivation of the formula for ordinary generating function for Ward-Horadam  sequence is a standard task and so we have (compare with  (5) in \cite[2009]{HorzumKocer}  by Tugba Horzum and
Emine G{\"{o}}kcen Kocer or put $r=1$  in Theorem 1 from  \cite[2003]{Pantelimon 2003}  by St$\breve{a}$nic$\breve{a}$) 

\begin{equation}
G_H(a,b;P,Q)[x] = \frac{a+(b-aP)x}{1-P \cdot x + Q \cdot x^2}=  \frac{a+(b-a[p+q])x}{1-P \cdot x + p\cdot q \cdot x^2}.   
\end{equation}
where from we decide an identification-abbreviation $$ G_H[x] \equiv  G_H(A,B;p,q)[x] \equiv G_H(a,b;P,Q)[x].$$

\vspace{0.2cm}

\noindent Naturally - in general  $H(A,B;p,q) \neq H(A,B;q,p)$. 

\vspace{0.2cm}
\noindent If $H(A,B;p,q) = H(A,B;q,p)$ we then call the  Ward-Horadam sequence $p,q$-symmetric and thus we arrive to Lucas  ($A=B=1$)
{\it Th\'eorie des Fonctions Num\'eriques Simplement Périodiques}  \cite[1878]{EdL}.

\vspace{0.2cm}

\noindent  In \cite[1878]{EdL} Edouard Lucas considers  Lucas sequence of the second kind   $V = \left\langle V_n\right\rangle_{n\geq 0}$ (second kind - see: \cite[1977, Part I]{KKK}) as well  as its till now well elaborated companion Lucas sequence of the first kind $ U = \left\langle U_n\right\rangle_{n\geq 0}$ (first  kind - see: \cite[1977, Part I]{KKK}) which gives rise in its turn to the  \textit{U}-binomials' recurrence (58) in \cite[1878]{EdL}. 

\noindent See then  \cite[1949]{JM},  \cite[1964]{TF},  \cite[1969]{Gould},  \cite[1989]{G-V}  or in  \cite[1989]{K-W}  etc. For example - for the relations  between  $U_{n+1}$  and $V_n$ number sequences and for the explicit form of ordinary generating functions of their powers see \cite[1977]{popov} by Blagoj  S. Popow, where the characteristic equation for both $H_n$ recurrent number sequences of second order is of general form $a\cdot z^2 + b\cdot z + c = 0$.

\vspace{0.2cm}

\noindent The  $p,q$-symmetric sequences from  \cite[1878]{EdL} i.e  the Lucas sequence of the second kind  ($A=B = 1$)

\begin{equation}
H_n(2,P;p,q)= V_n = p^n + q^n.   
\end{equation}
\noindent and  the Lucas sequence of the first kind  ($A= - B = 1$)

\begin{equation}
H_n(0,1;p,q)= U_n = \frac{p^n - q^n}{p-q},    
\end{equation}
where called by Lucas \cite[1878]{EdL} the \textit{simply periodic numerical functions} because of 

\vspace{0.1cm}

\noindent [quote] \textsl{at the start, the complete analogy of these symmetric functions with 
the circular and hyperbolic functions.} [end of quote].

\vspace{0.3cm}


\noindent \textbf{\textcolor{blue}{More ad Notation} 1.} 
\noindent The letters  a,b $a\neq b$ in \cite[1878]{EdL} denote the roots of the equation $x^2= P x -Q$  then  $(a,b)\mapsto (u,v)$ in \cite[2009]{Savage} and u,v stay there for the roots of the equation $x^2  =  \ell x -1$.

\noindent We shall use here the identification  $(a,b) \equiv (p,q)$  i.e.  $p,q$ denote the roots of $x^2= P x -Q$ as is common in "`Lucas $(p,q)$-people"' publications.


\noindent For Lucas $(p,q)$-people  then the following  $U$-identifications  are expediency natural:

\begin{defn}
     \begin{equation}
n_{p,q} = \sum_{j=0}^{n-1}{p^{n-j-1}q^j}=U_n =\frac{p^n - q^n} {p-q},\  0_{p,q}=U_0 = 0, \ 1_{p,q}=U_1 = 1,
    \end{equation}
\end{defn}
where  $p,q$  denote now the roots of the equation $x^2= P\cdot x - Q  \equiv  x^2= s x +t$ hence  $p+ q = s \equiv P$ , $pq = Q \equiv - t$ and the empty sum convention was used for $0_{p,q} = 0$.
\noindent Usually one assumes $p\neq q$. In general also $s\neq t$ - though  according to the context \cite[1989]{G-V}  $s=t$  may happen  to be the case of interest.

\vspace{0.2cm}

\noindent The Lucas $U$-binomial coefficients ${n \choose k}_U \equiv {n \choose k}_{p,q}$  are then defined as follows: (\cite[1878]{EdL},  \cite[1915]{Fon},  \cite[1936]{Ward},  \cite[1949]{JM},  \cite[1964]{TF}, \cite[1969]{Gould} etc.)

\vspace{0.3cm}

\begin{defn}
	Let  $U$ be as in \cite[1878]{EdL} i.e $U_n \equiv n_{p,q}$ then $U$-binomial coefficients for any $n,k \in \mathbb{N}\cup\{0\}$ are defined as follows
	\begin{equation}
		{n \choose k}_U \equiv {n \choose k}_{p,q} = \frac{n_{p,q}!}{k_{p,q}! \cdot (n-k)_{p,q}!} = \frac{n_{p,q}^{\underline{k}}}{k_{p,q}!}
	\end{equation}
	\noindent where $n_{p,q}! = n_{p,q}\cdot(n-1)_{p,q}\cdot ... \cdot 1_{p,q}$ and $n_{p,q}^{\underline{k}} = n_{p,q}\cdot(n-1)_{p,q}\cdot ...\cdot (n-k+1)_{p,q}$  and ${n \choose k}_U = 0$ for $k > n$.
\end{defn} 

\vspace{0.3cm}

\begin{defn}
Let $V$ be as in \cite[1878]{EdL} i.e $V_n = p^n +q^n$, hence $V_0 = 2$ and  $V_1 = p + q = s$. Then $V$-binomial coefficients for any $n,k \in \mathbb{N}\cup\{0\}$ are defined as follows
	\begin{equation}
		{n \choose k}_V =\frac{V_n!}{V_k!\cdot V_(n-k)!} = \frac{V_n^{\underline{k}}}{V_k!}
	\end{equation}
	\noindent where $V_n! = V_n \cdot V_{n-1}\cdot...\cdot V_1$ and $V_n^{\underline{k}}=V_n \cdot V_{n-1}\cdot ...\cdot V_{n-k+1}$ and ${n \choose k}_V = 0$ for $k > n$.
\end{defn}

\vspace{0.2cm}

\noindent One automatically generalizes number $F$-binomial coefficients' array to functions $F(x)$-\textbf{multinomial} coefficients' array (see  \cite[2004]{akk3 2004} and references to umbral calculus therein) while for \textit{number sequences} $F = F(x=1)$ the $F$-multinomial coefficients see p. 40 in \cite[1913]{Carmichel} by  Robert Daniel Carmichel , see  \cite[1936]{Ward} by Morgan Ward and \cite[1969]{Gould}  by Henri W. Gould.

\vspace{0.2cm}

\noindent In  \cite[1979]{Shannon 1979 multi})   Anthony G. Shannon considers the special case of number sequences, where the  $F_n = H^{(r)}_{s,n},  n =0,1,2,...$, $1 \leq s \leq r$, which are constituting $r$ basic number sequences - satisfy a linear homogeneous recurrence relation of the order $r$:
 
 $$ H^{(r)}_{s,n} = \sum_{j=1}^r\left(-1\right)^{j+1}P_{r,j}\cdot H^{(r)}_{s,n-j}, \ n>r,  $$
 with initial conditions $H^{(r)}_{s,n}= \delta_{s,n}$ for  $1 \leq n \leq r$ and where $P_{r,j}$ are arbitrary integers.

\noindent For $r=2$ the above coincides  with this note recurrence (1)  where  $H^{(2)}_{1,1}= H_1 = 1$, $H^{(2)}_{1,2}= H_0 = 0$  for  $s=1$  and  $H^{(2)}_{2,1}= H_0 = 0$, $H^{(2)}_{2,2}= H_1 = 1$  for  $s=2$.  Because of that latter Anthony G. Shannon designates the number sequence $\left\langle H^{(r)}_{r,n} \right\rangle \equiv \left\langle u^{(r)}_n \right\rangle $ as the fundamental sequence by the analogy with Lucas' second-order fundamental sequence  $\left\langle H^{(2)}_{r,n} \equiv U_n \right\rangle$.
 
\vspace{0.1cm}

\noindent Anthony G. Shannon derives then in \cite[1979]{Shannon 1979 multi}  the recursion for $u^{(r)}$-multinomial coefficients (see Theorem on page 346), where there specifically  for $r=2$ the number sequence  $u^{(2)}$ coincides with  the fundamental Lucas sequence (i.e. the Lucas sequence of the first kind).

\vspace{0.1cm}

\noindent As for other here relevant works see \cite[1991]{SatoShi} by Shiro Ando and  Daihachiro Sato.  Note also  $x$-Fibonomial coefficients array from  \cite[2001]{Rich 2001}  by  Thomas M. Richardson.

\vspace{0.1cm}

\noindent \textbf{$\psi(x)$-multinomial Remark.}

\vspace{0.1cm}

\noindent Considerations in \cite[1979]{Shannon 1979 multi}  by  Anthony G. Shannon  and an application in \cite[2001]{Rich 2001} by Thomas M. Richardson as well as relevance to  umbral calculus 
\cite[2003]{akk ITSF 2003} or  \cite[2001]{akk ITSF 2001} constitute  motivating circumstances for considering now \textit{functions'} $F(x)$-\textbf{binomial}  and  $F(x)$-\textbf{multinomial} coefficients. This is the case ($F(x)=\psi(x)$) for example in  \cite[2004]{akk3 2004}  \textbf{wherein we read:} 

$$
\left(x_{1}+_{\psi}x_{2}+_{\psi} \ldots +_{\psi}x_{k}\right)^{n}=
\sum_{\begin{array}{l} s_{1},\ldots s_{k}=0\\
s_{1}+s_{2}+\ldots +s_{k}=n
\end{array}}^{n}\binom{n}{s_{1},\ldots ,s_{k}}_{\psi}x_{1}^{s_{1}}\ldots x_{k}^{s_{k}}
$$
where
$$\binom{n}{s_{1},\ldots ,s_{k}}_{\psi}=\frac{n_{\psi}!}{(s_{1})_{\psi}!\ldots (s_{k})_{\psi}!}.$$
Here above   the $u$-{\em multinomial} number sequence formula from \cite[1936]{Ward} by Morgan Ward is extended mnemonically to $\psi(x)$ function sequence definition of shifting $x$ arguments formula written in Kwa\'sniewski upside-down notation (see for example \cite[2003]{akk ITSF 2003},  \cite[2001]{akk ITSF 2001}, \cite[2002]{akk Lodz2002},  \cite[2005]{akk Allahabad  2005} for more ad this notation).    

\noindent  The above formulas introduced  in case of number sequences in \cite[1936]{Ward} by Morgan Ward were recalled from \cite[1936]{Ward}  by Alwyn F. Horadam and  Anthony G.  Shannon  in \cite[1976]{HorSha}  were the Ward Calculus of sequences framework (including umbral derivative and corresponding exponent) was used to enunciate two types of  Ward's Staudt-Clausen theorems pertinent to this general calculus of number sequences.

\noindent The end of \textbf{$\psi(x)$-multinomial Remark.}

\vspace{0.2cm}

\noindent Recently polynomial-fibonomial coefficients in  \cite[2008] {Cigler 2008} by Johann Cigler  (i.e. $x$-Fibonomial coefficients  from  \cite[2001]{Rich 2001}) appear naturally in derivation of recurrence relations for powers of $q$-Fibonacci polynomials  by Johann Cigler. Let us then come over to this multinomiality closing definition.

\begin{defn}
Let $F(x)$ be any natural, or complex  numbers' non zero valued functions' sequence i.e. $F_n(x)\in\mathbb{N}$ or and $F_n(x)\in\mathbb{C}$. The \textbf{$F(x)$-multinomial coefficient} is then identified with the symbol

\begin{equation}
	{n \choose {k_1,k_2,...,k_s}}_{F(x)} = \frac{F_n(x)!}{F_{k_1}(x)!\cdot ... \cdot F_{k_s}(x)!}
\end{equation}

\vspace{0.2cm}
\noindent where $k_i\in\mathbb{N}$ and $\sum_{i=1}^{s}{k_i} = n$ for $i=1,2,...,s$. Otherwise it is equal to zero,
\end{defn}
\noindent and where $F_r(x)! = F_r(x) \cdot F_{r-1}(x)\cdot...\cdot F_1(x)$.

\vspace{0.2cm}

\noindent Naturally  for any natural $n,k$ and $k_1+...+k_m=n-k$ the following holds

\begin{equation} 
{n \choose k}_{F(x)} \cdot {n-k \choose {k_1,k_2,...,k_m}}_{F(x)}  = {n \choose {k,k_1,k_2,...,k_m}}_{F(x)}, 
\end{equation}
 
$${n \choose {k_1,k_2,...,k_m}}_{F(x)}={n \choose k_1}_{F(x)}{n-k_1\choose{k_2}}_{F(x)}\cdots {n-k_1- \cdots -k_{m-1}\choose{k_m}}_{F(x)}.$$


\vspace{0.3cm}


\noindent \textbf{\textcolor{blue}{More ad Notation} 2.} 

Does notation of items' representing - matter?

\vspace{0.2cm}

\noindent As a \textcolor{blue}{\textbf{Motto}} for an experienced answer we propose:

\vspace{0.1cm}

\noindent \textsl{Science is a language}  (from \cite[1996]{Doron ZeilBerger 1996}   by  Doron Zeilberger); 

\vspace{0.1cm}

\noindent \textit{Mathematical notation evolves like all languages}  (from \cite[1992]{Knuth 1992} by Donald Ervin Knuth).

\vspace{0.2cm}

\noindent We shall use further on the traditional , XIX-th century rooted notation under presentation in spite of being inclined to quite younger notation from \cite[2009]{BSCS}  by  Bruce E. Sagan and Carla  D.  Savage.  This wise, economic notation is ready for straightforward 	record of combinatorial interpretations and combinatorial interpretations' substantiation in terms of popular text book tiling model since long ago used for example to visualize recurrence for Fibonacci-like sequences ; see for example \cite[1989]{GKP 1989}  by Ronald Graham, Donald Ervin  Knuth, and Oren Patashnik.  
The translation from Fran\c{c}ois  \'Edouard Anatole Lucas via Dov Jarden and Theodor Motzkin  notation \cite[1949]{JM} and notation  of Bruce E. Sagan and Carla  D.  Savage \cite[2009]{BSCS} is based 
on the succeeding identifications: the symbol used for $U$-\textbf{binomials} is  $\left\{... \right\}$  in place of  $\left(...\right)_U$. 

\noindent The would be symbol  for $V$-\textbf{binomials} i.e. $P=\left\langle ... \right\rangle$  in place of  $\left(...\right)_V$  is not considered at all in  \cite{BSCS} while  

$$ \left\{n \right\} \equiv U_n \equiv n_{p,q}  , \ \left\langle n \right\rangle \equiv V_n.$$
In Bruce E. Sagan and Carla  D.  Savage  notation  we would then write down the  fundamental and primordial sequences' binomial coefficients  as follows.

\vspace{0.3cm}

\begin{defn}
	Let  $\left\{n\right\}$ be fundamental Lucas sequence as in \cite[1878]{EdL} i.e $\left\{n\right\}\equiv U_n \equiv n_{p,q}$ then $\left\{n\right\}$-binomial coefficients for any $n,k \in \mathbb{N}\cup\{0\}$ are defined as follows
	\begin{equation}
		F\left\{n,k \right\} = \left\{ n \atop k\right\}_{p,q} = \frac{\left\{n\right\}!}{\left\{k\right\}! \cdot \left\{n-k\right\}!} = \frac{\left\{n\right\}^{\underline{k}}}{\left\{k\right\}!}
\end{equation}
where $\left\{n \right\}! = \left\{n \right\}\cdot\left\{n-1\right\}\cdot...\cdot \left\{1\right\}$  
and $\left\{n \right\}^{\underline{k}} = \left\{n\right\}\cdot\left\{n-1\right\}\cdot ...\cdot\left\{n-k+1\right\}$ and $F\left\{n,k \right\} = 0$ for $k > n$.
\end{defn}

\begin{defn}
	Let  $\left\langle n \right\rangle $ be primordial Lucas sequence as in \cite[1878]{EdL} i.e $\left\langle n \right\rangle \equiv V_n $ then $ \left\langle n \right\rangle $-binomial coefficients for any $n,k \in \mathbb{N}\cup\{0\}$ are defined as follows
	\begin{equation}
P\left\langle n,k \right\rangle = \left\langle  n \atop k\right\rangle_{p,q} =\frac{\left\langle n \right\rangle!}{\left\langle k \right\rangle! \cdot \left\langle n-k \right\rangle!} = \frac{\left\langle n \right\rangle^{\underline{k}}}{\left\langle k \right\rangle!},
\end{equation}
where $\left\langle n \right\rangle! = \left\{n \right\}\cdot\left\{n-1\right\}\cdot...\cdot \left\{1\right\}$ and $P\left\langle n,k \right\rangle = 0$ for $k > n$. 
and $\left\{n \right\}^{\underline{k}} = \left\{n\right\}\cdot\left\{n-1\right\}\cdot ...\cdot\left\{n-k+1\right\}.$
\end{defn} 
The above consequent symbols $\left\{ n \atop k\right\}_{p,q}$  and  $\left\langle  n \atop k\right\rangle_{p,q}$     are occasionally  - in  not exceptional  conflict - with second kind Stirling numbers notation and Euler numbers notation respectively in the spirit of \cite[1989]{GKP 1989} what extends on both  $p,q$ - extensions' notation.  

\noindent Regarding  the symbol $\left\{ n \atop k\right\}_{p,q}$ one draws the attention of a reader to  \cite[1967]{Hoggatt 1967} where Verner Emil Hoggatt Jr.    considers the $U$-binomial coefficients ${n \choose k}_U \equiv {n \choose k}_{p,q}$  denoting them as  $ \left\{ n \atop k\right\}$ where $\left\{u_n\right\}_{n\geq 0}= U $  with  $U$ being the fundamental Lucas sequence. 
The author of \cite[1967]{Hoggatt 1967} derives among others also recurrences  [see  this note formula (43)]

$$\left\{ n \atop k\right\} = u_{k+1}\cdot \left\{ n-1 \atop k\right\} + t\cdot u_{n-k-1}\cdot \left\{ n-1 \atop k-1\right\} ,$$ 

$$\left\{ n \atop k\right\} =  t\cdot u_{k-1}\cdot \left\{ n-1 \atop k\right\} + u_{n-k+1}\cdot \left\{ n-1 \atop k-1\right\} ,$$ 

$$\left\{ n \atop k\right\}_{F(x)} = f(x)_{k+1}\cdot \left\{ n-1 \atop k\right\}_{f(x)} + 1 \cdot f(x)_{n-k-1}\cdot \left\{ n-1 \atop k-1\right\}_{F(x)} ,$$ 
where  $ F(x)= \left\langle f(x)_n \right\rangle_{n>0}$ denotes the sequence of Fibonacci polynomials.  

\vspace{0.2cm}

\noindent These recurrences are to be placed side by side  with corresponding recurrence formulas in $\left\{...\right\}$ notation as for example the recurrence of  the Theorem p.346  for $r=2$ in \cite[1979]{Shannon 1979 multi}) by  Anthony G. Shannon  or the  recurrence (10.3) in  \cite[1989]{G-V} by  Ira  M. Gessel  and Xavier G\'erard Viennot or with the Proposition 2.1  in \cite[2009]{BSCS}  by Bruce E. Sagan and Carla  D.  Savage.  

\vspace{0.2cm}

\noindent Anthony G. Shannon designates the recurrence of  the Theorem p.346 in \cite[1979]{Shannon 1979 multi}) as the recurrence relation for multinomial coefficients

$$ 	{n \choose {s_1,s_2,...,s_r}}_u = \frac{u^{(r)}_n!}{u^{(r)}_{s_1}!\cdot ... \cdot u^{(r)}_{s_r}!} ,$$
where   $u^{(r)}_n \equiv U^{(r)}_{r,n+r}$  and naturally $\sum_{k=1}^r s_k = n $.

\vspace{0.1cm}

\noindent For  $ \left\{ n \atop k\right\}$ corresponding notation see then  also\textbf{:} \cite[1969]{Gould} by Henri W.Gould, \cite[1976]{HorSha} by Alwyn F. Horadam and  Anthony G.  Shannon,  \cite[1989]{G-V} by Ira  M. Gessel and Xavier G\'erard Viennot , \cite[2005]{ST 2005}, \cite[2005]{ST 2005 2}, \cite[2006]{ST 2006}  by Jaroslav Seibert  and Pavel Trojovsk\'y and \cite[2007]{T 2007}  by  Pavel Trojovsk\'y.

\vspace{0.2cm}

\noindent Whereas as in the\textit{ subset-subspace problem} ( Example [Ex. q* ; 6] in subsection 4.3.)  we rather need another natural notation. Namely  for  $q \neq 0$  introduce  $q* = \frac{p}{q}$ and observe that 

$$ {n \choose k}_U \equiv {n \choose k}_{p,q} = q^{k(n-k)} {n \choose k}_{1,q*} \stackrel{q*\mapsto 1}{\rightarrow}   \fnomialF{n}{k}. $$    
The $V$-binomial  $P\left\langle n,k \right\rangle = \left\langle  n \atop k\right\rangle_{p,q} \equiv \fnomialF{n}{k}{V}$ is not considered in \cite[2010]{BSCS} - neither recurrences for $\fnomialF{n}{k}{V}$ coefficients  are derived in earlier publications - up to knowledge of the present author .




\section{$H(x)$-binomial coefficients' recurrence}

\vspace{0.2cm}

\noindent \textbf{3.1.}  Let us recall convention resulting from (3).

\vspace{0.2cm}
\noindent Recall. The general solution of (3):   $$ H(x) \equiv H(a(x),b(x);s(x),t(x)) = \left\langle H_n(x)\right\rangle _{n\geq 0}$$ 
is being called  throughout this paper - \textbf{Ward-Horadam functions' sequence}.

\vspace{0.2cm}

\noindent  From the characteristic equation of the recurrence (3)

\begin{equation}
z^2 - s(x) \cdot z  - t(x) =0
\end{equation}
we readily see that  for $ H_0 = a(x), \:H_1 = b(x)$  and $ \;n\geq 0 ,$

\begin{equation}
H_n(x) \equiv H_n(a(x),b(x);p(x),q(x)) = A(x)p(x)^n + B(x)q(x)^n,                                                                                     
\end{equation}
where  $p(x),q(x)$ are roots of (20) and we have assumed that  $p(x)\neq q(x)$  as well as that  $\frac{p(x)}{q(x)}$ are not roots of unity  (see  for example Lemma 1 in \cite[1980]{Beukers} 
by Fritz Beukers). Naturally:
 
\begin{equation}
A(x) = \frac{b(x)-q(x)a(x)}{p(x)-q(x)} \ , \;    B = - \frac{b(x)-p(x)a(x)}{p(x)-q(x)}.   
\end{equation}
\noindent hence we may and we shall use the following  conventional  identifications-abbreviations

\begin{equation}
H(x) \equiv H(a(x),b(x);s(x),t(x)) \equiv  H(A(x),B(x);s(x),t(x)).
\end{equation}

\vspace{0.2cm}

\noindent As for the case $p(x)=q(x)$  included see for example Proposition 2.7 in  \cite[2009]{He-Siue}.

\vspace{0.2cm}

\noindent Another  explicit formula for Ward-Horadam functions sequences is the mnemonically  extended formula (9) from  \cite[2009]{HorzumKocer} by Tugba Horzum and Emine G{\"{o}}kcen  Kocer,  where here down we use contractually the following abbreviations:

\vspace{0.2cm}

\noindent $H_n(x) \equiv H_n(a(x),b(x);s(x),t(x))$ , $a(x) \equiv a $, $b(x) \equiv b $, $s(x) \equiv s $ and  $t(x) \equiv t $

\begin{equation}
H_n(x)= a\sum_{0\leq k \leq \left\lfloor \frac{n}{2}\right\rfloor} \binom{n-k}{k}s^{n-2k}t^k + \left( \frac{b}{s}-a \right)\sum_{0\leq k \leq \left\lfloor \frac{n-1}{2}\right\rfloor} \binom{n-k-1}{k}s^{n-2k}t^k                                                                          
\end{equation}

\vspace{0.2cm}

\noindent \textbf{Note and compare}. The recurrence  (1.1) and (1.2) in \cite[1996]{1 Horadam 1996} by Alwyn  F. Horadam defines a polynomials' subclass of Ward-Horadam functions sequences defined by (3).
The standard Jacques Binet  form (1.8) in \cite[1996]{1 Horadam 1996} of the  recurrence  (1.1) and (1.2) solution for Ward-Horadam polynomials  sequences in \cite[1996]{1 Horadam 1996} is the standard Jacques Binet  form (19), (20) of the  recurrence  (3) solution for Ward-Horadam functions sequences.

\vspace{0.1cm}

\noindent The recurrence  (2.23) in \cite [2009]{He-Siue}  by  Tian-Xiao He and   Peter Jau-Shyong Shiue defines exactly the class of Ward-Horadam functions' second order sequences 
and this paper standard Jacques Binet  form (21), (22) of the  recurrence  (3) solution for Ward-Horadam functions sequences $H(x)$ constitutes the content of their Proposition 2.7. - as has been  mentioned earlier.  No recurrences for $H(x)$-binomials neither for $H(x=1)$-binomials are considered.

\vspace{0.2cm}

\noindent \textbf{On Binet Formula - Historical Remark.}  We just quote Radoslav Rasko Jovanovic's information from   $$ http://milan.milanovic.org/math/english/relations/relation1.html :$$

\begin{quot} Binet's Fibonacci Number Formula was derived by Binet in 1843 although the result was known to Euler and to Daniel Bernoulli more than a century ago. ... It is interesting that A de Moivre (1667-1754) had written about Binet`s Formula, in 1730, and had indeed found a method for finding formula for any general series of numbers formed in a similar way to the Fibonacci series.
\end{quot} 
See also the book  \cite[1989]{Vajda 1989}    by  Steven Vajda. 

\vspace{0.2cm}

\noindent \textbf{3.2.}  The authors of \cite{Savage} provide an easy proof of  a general observation named there deservedly Theorem (Theorem 17) which extends automatically to the statement that  the following  recurrence holds for the general case of ${r+s \choose r,s}_{H(x)}$  $H(x)$-binomial array in multinomial notation.

\vspace{0.2cm}

\begin{theoremn}
Let us admit shortly the abbreviations: $g_k(r,s)(x) = g_k(r,s)$ , $k=1,2$.  Let $s,r>0$. Let $F(x)$ be any zero characteristic field nonzero valued functions' sequence ($F_n(x) \neq 0$). Then

\begin{equation}
{r+s \choose r,s}_{F(x)}=g_1(r,s)\cdot{r+s-1 \choose r-1,s}_{F(x)}+g_2(r,s)\cdot{r+s-1 \choose r,s-1}_{F(x)}
\end{equation}
where   $ {r \choose r,0}_{F(x)} = {s \choose 0,s}_{F(x)} =1$ and 
\begin{equation}
 F(x)_{r+s} =  g_1(r,s) \cdot F(x)_r   +  g_2(r,s) \cdot F(x)_s.                               
\end{equation}
are equivalent.
\end{theoremn}

\vspace{0.3cm}

\noindent \textbf{On the way historical note}
\noindent  Donald  Ervin  Knuth   and    Herbert Saul Wilf in \cite[1989]{K-W} stated that  Fibonomial coefficients and the recurrent relations for them appeared already in 1878  Lucas work (see: formula (58)  in \cite[1878]{EdL}  p. 27 ;  for $U$-binomials which "Fibonomials" are special case of). More over on this very p. 27 Lucas formulated a conclusion from his (58) formula which may be  stated in notation of this paper formula (2) as follows: \textit{if}  $s,t  \in\mathbb{Z}$ and $H_0 =0$ , $H_1 = 1$  \textit{then} $H \equiv U$  and $ \fnomialF{n}{k}{U} \equiv \fnomialF{n}{k}{n_{p,q}} \in\mathbb{Z}$.  

\vspace{0.1cm}

\noindent Consult for that also the next century references: \cite[1910]{Bachmann 1910}  by  Paul Gustav Heinrich Bachmann  or later on - \cite[1913]{Carmichel}  by  Robert Daniel Carmichel  [p. 40]  or  \cite[1949]{JM} by Dov Jarden and Theodor Motzkin where in all quoted positions it was also shown that $n_{p,q}$ - binomial coefficients are integers  - for $p$ and $q$ representing  distinct roots of this note characteristic equation (5).

\vspace{0.2cm}

\noindent Let us take an advantage to note that  Lucas  Th\'eorie des Fonctions Num\'eriques Simplement P\'eriodiques i.e. investigation of exactly  fundamental $U$  and  primordial $V$ sequences constitutes  the far more non-accidental context  for binomial-type coefficients exhibiting their relevance at the same time to number theory and  to hyperbolic trigonometry  (in addition to \cite[1878]{EdL} see for example \cite{akk1999},  \cite{ Bakk} and  \cite{Bajguz}).  

\vspace{0.1cm}

\noindent It seems to be the right place now to underline that the  \textit{addition formulas} for Lucas sequences below with respective hyperbolic trigonometry formulas and also consequently $U$-binomials'recurrence formulas - stem from commutative  ring $R$ identity: $(x-y)\cdot (x+y) \equiv x^2 - y^2, x,y \in R$.

\vspace{0.3cm}

\noindent Indeed. Recall the characteristic equation notation:  $z^2 =\textcolor{blue}{z^2 = s \cdot x + t}$ as is common in many publications with the restriction:  $p,q$ roots are distinct. Recall that $p+q=s$  and $p\cdot q=-t$. Let $\Delta = s^2 - 4t$. Then $\Delta = (p-q)^2$. Hence we have as in \cite[1878]{EdL}

\begin{equation}
2\cdot U_{r+s} =  U_r V_s  +  U_s V_r ,\ \ \ \ 
2\cdot V_{r+s} =  V_r V_s  +  \Delta \cdot U_s U_r. 
\end{equation}
Taking here into account the  $\textbf{U}$-\textcolor{blue}{addition formula} i.e. the first of two trigonometric-like $L$-addition formulas (42) from \cite[1878]{EdL} ($L = U,V$)  
one readily recognizes that  the $U$-binomial recurrence from the Corollary 18 in  \cite[2009]{Savage} is  the $U$-binomial recurrence (58) \cite[1878]{EdL}  which may be rewritten after Fran\c{c}ois  \'Edouard Anatole Lucas in multinomial notation and stated as follows: \textit{according  to the Theorem 2a  below the following is true}:

$$2\cdot U_{r+s} =  U_r V_s  +  U_s V_r$$ 
\textit{is equivalent to}

\begin{equation}
2\cdot{r+s \choose r,s}_{n_{p,q}}= V_s \cdot{r+s-1 \choose r-1,s}_{n_{p,q}}+V_r \cdot{r+s-1 \choose r,s-1}_{n_{p,q}}.
\end{equation}
\noindent  To this end see also Proposition 2.2. in \cite[2009]{BSCS} and compare it with both (28) and  Example 3. below.

\vspace{0.2cm}

\noindent However there is no companion  $V$-binomial recurrence i.e. for ${r+s \choose r,s}_V$ neither in  \cite[1878]{EdL}  nor in \cite[2009]{Savage} as well as all other quoted papers - up to knowledge of this note author.  

\vspace{0.2cm}

\noindent Consequently  then there is no overall $H(x)$-binomial recurrence neither in  \cite[1878]{EdL}  nor in \cite{Savage} (2009) as well as all other quoted papers  except for Final remark : p.5 in \cite[2010]{akk2010} up to  knowledge of the present author.  

\vspace{0.3cm}

\noindent \textbf{The End} of \textit{the on the way historical note}.

\vspace{0.4cm}

\noindent The looked for  $H(x)$-binomial recurrence  (29) accompanied by (30) might be  then given right now in the form of (25) adapted to  - Ward-Lucas functions'sequence case notation while keeping in mind that of course the expressions for $h_k(r,s)(x)$,  $k=1,2$ below are designated by this $ F(x) = H(x)$ choice and as a matter of fact are appointed by the recurrence (3). 

\vspace{0.3cm}

\noindent For the sake of commodity to write down an observation to be next (named deservedly ? Theorem) let us admit shortly the abbreviations: $h_k(r,s)(x) = h_k(r,s)= h_k$ , $k=1,2$. Then for $H(x)$ of the form (21) we evidently have what follows.

\vspace{0.6cm}

\noindent \textbf{\textcolor{blue}{Theorem 2a}.}

\begin{equation}
{r+s \choose r,s}_{H(x)} = h_1(r,s){r+s-1 \choose r-1,s}_{H(x)} + h_2(r,s){r+s-1 \choose r,s-1}_{H(x)},
\end{equation}
where  ${r \choose r,0}_{H(x)} = {s \choose 0,s}_{H(x)} =1,$  \textit{is equivalent to}

\begin{equation}
 H_{r+s}(x) =  h_1(r,s)H_r(x)   +  h_2(r,s) H_s(x).                               
\end{equation}
where  $H_n(x)$ is explicitly given by (21)  and (22). 
\vspace{0.1cm}
\noindent \textbf{The end} of the Theorem 2a.

\vspace{0.2cm}

\noindent There might be various  $h_1(r,s)(x)=h_1$ and $ h_2(r,s)(x)_s =h_2 $ solutions of  (30) and  (21). Compare (38)  in Example 1 with (42) in Example 3  below. 

\vspace{0.1cm}
\noindent As the possible  $h_1(r,s)(x)=h_1$ and $ h_2(r,s)(x)_s =h_2 $ formal solutions of  (30) and  (21)  we just may take

\begin{equation}
h_1(r,s)(x)= \frac{A(x)\cdot p(x)^{r+s}} {A(x) \cdot p^r + B(x) \cdot q(x)^r},   
\qquad
h_2(r,r)(x) = \frac{B(x) \cdot q(x)^{r+s}} {A(x) \cdot p^s + B(x) \cdot q(x)^s}.
\end{equation}
As  another possible $h_1(r,s)(x)$ and $ h_2(r,s)(x)_s =h_2 $ solutions of  (30) and  (21) we may take:   \textbf{for $r \neq s$}

\begin{equation}
h_1(r,s) \cdot (p((x)^r q(x)^s  -  q(x)^r p(x)^s ) = p(x)^{r+s} q(x)^s  -  q(x)^{r+s} p(x)^s,
\end{equation}

\begin{equation}
h_2(r,s) \cdot (q(x)^r p(x)^s  -  p(x)^r q(x)^s ) = p(x)^{r+s} q(x)^r  -  q(x)^{r+s} p(x)^r. 
\end{equation}
\textit{while for }$r=s$  apply formula (31)  with $r=s$.

\vspace{0.2cm}

\noindent Usually the specific features of  particular cases of (21)  and (22) allow one to infer the particular form of (30) hence the form of $h_1(r,s)(x)=h_1$ and $ h_2(r,s)(x)_s =h_2 $.

\vspace{0.2cm}

\noindent As we soon shall see below (Theorem 2b), there is a specific way out from being  bothered by this not illusory obstacle in $V$-binomials case. In general case one may always use the  Georges Fonten\'e  recurrence \cite[1915]{Fon}  which  however seems not to have combinatorial interpretation except for natural and $0$  numbers $n$ and their  $q$-extensions $n_q = n_{q,q}$ ; $n \geq 0$. See Example 2 below.

\vspace{0.3cm}

\noindent \textbf{3.3. Three special cases examples.}

\vspace{0.2cm}

\noindent \textbf{Example 1.} This is a particular case of the Theorem 2a. 
\vspace{0.1cm}

\noindent The recurrent relations (13) and (14) in  Theorem 1 from  \cite[2008]{Corsino} by Roberto Bagsarsa Corcino  for $n_{p,q}$-binomial coefficients are special cases of this paper formula (29) as well as  of Th. 17 in \cite{Savage} with straightforward identifications of $g_1, g_2$  in (13)  and  in (14) in \cite{Corsino} or in this paper recurrence (30) for $H(x=1) = U[p,q]_n = n_{p,q}$ sequence. Namely, recall here now in multinomial notation this  Theorem 1 from  \cite[2008]{Corsino}   by Roberto Bagsarsa Corcino:

\begin{equation}
{r+s \choose r,s}_{p,q} = q^r{r+s-1 \choose r-1,s}_{p,q} + p^s{r+s-1 \choose r,s-1}_{p,q},
\end{equation}

\begin{equation}
{r+s \choose r,s}_{p,q} = p^r{r+s-1 \choose r-1,s}_{p,q} + q^s{r+s-1 \choose r,s-1}_{p,q},
\end{equation}
which is equivalent to

\begin{equation}
(s + r)_{p,q} = p^s r_{p,q} +  q^r s_{p,q} = (r + s)_{q,p} = p^r s_{p,q} +  q^s  r_{p,q},
\end{equation}
what  might be at once seen proved by noticing that $$p^{r+s} - q^{r+s} \equiv p^s \cdot(p^r - q^r) +   q^r \cdot(p^s - q^s).$$  Hence those mentioned straightforward identifications follow:

\begin{equation}
 g_1 = q^r ,\ \  g_2 =  p^s \;or\ \  g_1 = p^r ,\ \  g_2 =  q^s. 
\end{equation}
The recurrence (36) in Lucas notation reads
\begin{equation}
U_{s+r} = p^s U_r +  q^r U_s = U_{r+s} = p^r U_s +  q^s U_r. 
\end{equation}
Compare it with equivalent recurrence (42) from the \textbf{Example 3.}  in order to notice that  both $h_1$  and $h_2$ functions are different from case to case  of recurrence (30) \textbf{equivalent realizations}.

\vspace{0.2cm}

\noindent  Compare this example based on Theorem 1  in \cite[2008]{Corsino} by Roberto Corcino with  with \cite[2008]{MD2} \textbf{v[1]} by Maciej Dziemia\'nczuk (see there (1) and (2) formulas), and track as well  -  the simple combinatorial proof of the Corollary 3   in  \cite[2009]{MD2} \textbf{v[2]})  by Maciej Dziemia\'nczuk.

\vspace{0.2cm}

\noindent \textbf{Example 2.} This is a particular case of the Theorem 1. 
\vspace{0.1cm}

\noindent Now let $A$ be any natural numbers' or even complex numbers' valued sequence. One readily sees that also (1915 year from) Fonten\'e recurrence for Fonten\'e-Ward generalized $A$-binomial coefficients i.e. equivalent identities (6) , (7) in \cite{Gould} \textbf{are special cases of} this paper formula \textbf{(25)} as well as  of Th. 17 in \cite{Savage} with straightforward identifications of $h_1, h_2$  in this paper formula (25)  while this paper recurrence  (26) becomes trivial identity.

\noindent Namely, the identities (6) and (7) from \cite[1969]{Gould} read correspondingly:

\begin{equation}
{r+s \choose r,s}_A = 1 \cdot {r+s-1 \choose r-1,s}_A + \frac{A_{r+s} - A_r}{A_s}{r+s-1 \choose r,s-1}_A,
\end{equation}

\begin{equation}
{r+s \choose r,s}_A = \frac{A_{r+s} - A_s} {A_{r}} \cdot {r+s-1 \choose r-1,s}_A +  1 \cdot {r+s-1 \choose r,s-1}_A,
\end{equation}
where $p \neq q$ and  $ {r \choose r,0}_L = {s \choose 0,s}_L =1.$ And finally we have tautology identity

\begin{equation}
A_{s + r} \equiv \frac{A_{r+s} - A_s}{A_r}\cdot A_r +  1 \cdot A_s.
\end{equation}
Example 2. becomes the general case of the Theorem 1. if we allow $A$ to represent any zero characteristic field nonzero valued functions' sequence: $A = A(x) = \left\langle A_n(x)\right\rangle_{n \geq 0}, \;  A_n(x) \neq 0$).   In particular we may put $A_n(x) = H_n(x)\equiv H_n(s(x),t(x);A(x),B(x))$  to put into play also $V$-binomials; $V_n \equiv H_n(s,t;1,-1).$

\vspace{0.2cm}

\noindent \textbf{Example 3.} This is a particular case of the Theorem 2a. 
\vspace{0.1cm}

\noindent The first example above is cognate to this third example in apparent way  as might readily seen from Fran\c{c}ois  \'Edouard Anatole Lucas papers \cite[1878]{EdL} or more  recent article  \cite[2001]{Hu-Sun}  by Hong Hu and Zhi-Wei Sun  ; (see also $t=s$ case in \cite[1989]{G-V} by Ira  M. Gessel and Xavier G\'erard Viennot on pp.23,24 .) In order to experience  this let us start to consider now the number $H(x=1)=U$ Lucas fundamental sequence \textit{fulfilling} this note recurrence  (\textbf{2}) with  $U_0 = 0$  and $U_1 = 1$  as  introduced in \cite[1878]{EdL} and then - for example  considered in \cite[2001]{Hu-Sun}. There in \cite[2001]{Hu-Sun} by Hong-Hu and Shi-Wei Sun - as a matter of fact - a kind of "pre-Theorem 17"  from \cite[2009]{Savage} is latent in the  proof of  Lemma 1 in  \cite{Hu-Sun}.  We rewrite  Lemma 1  by Hong-Hu and Shi-Wei Sun  in multinomial notation and an arrangement convenient for our purpose here using sometimes abbreviation $U_n(p,q)\equiv U_n$. 

\noindent (Note that the  \textit{addition formulas} for Lucas sequences hence consequently $U$-binomials'recurrence formulas  \cite[1878]{EdL} as well as  $(p-q)\cdot(p^{j+k} - q^{j+k}) \equiv (p^{k+1} - q^{k+1})\cdot(p^j - q^j) - p\cdot q (p^{j-1} - q^{j-1} \cdot(p^k - q^k)$ - stem from commutative ring $R$ identity: $(x-y)\cdot (x+y) \equiv x^2 - y^2, x,y \in R$.)

\vspace{0.2cm}

\noindent And so for   $p \neq q$ and bearing in mind that  $p\cdot q = - t$ - the following is true: 

\vspace{0.2cm}

\noindent \textit{The identity }(42) \textit{equivalent to }
$$(p-q)\cdot(p^{j+k} - q^{j+k}) \equiv (p^{k+1} - q^{k+1})\cdot(p^j - q^j) - p\cdot q (p^{j-1} - q^{j-1}) \cdot(p^k - q^k)$$  

\begin{equation}
U_{j + k}(p,q) = U_{k+1} \cdot U_j(p,q) + t \cdot U_{j-1}\cdot U_k(p,q)
\end{equation}
\noindent\textit{is equivalent to}

\begin{equation}
{j+k \choose j,k}_U = U_{k+1} \cdot {j+k-1 \choose j-1,k}_U + t \cdot U_{j-1}\cdot{j+k-1 \choose j,k-1}_U,
\end{equation}
\textit{where} $p,q$  \textit{are the roots of} (5) \textit{and correspondingly the above Lucas fundamental sequence} $H_n = U_n(p,q)$ i.e.  $U_0 = 0$  and $U_1 = 1$ \textit{is given by its Binet form} (6),(7). 

\vspace{0.2cm}

\noindent\textbf{Compare} now (42) with equivalent this note recurrence (38) in order to notice that  both $h_1$  and $h_2$ functions are different from case to case  of recurrence (30) \textit{equivalent realizations}.  

\vspace{0.1cm}

\noindent \textbf{Compare} then this paper recurrence formula (42)  with recurrence formula (4) in \cite[2009]{BSCS} or the  recurrence (F) from \cite[1967]{Hoggatt 1967} by Verner Emil Hoggatt, Jr. (see the end of this note Section 2) or $r=2$ case of the the recursion for $u^{(r)}$-multinomial coefficients (see Theorem on page 346)   in \cite[1979]{Shannon 1979 multi}  by  Anthony G. Shannon (see this note Section 2).   

\vspace{0.1cm}

\noindent \textbf{Compare} this paper recurrence formula (43)  with Proposition 2.2.  in \cite[2009]{BSCS}  by  Bruce E. Sagan and  Carla  D.  Savage. 

\vspace{0.1cm}

\noindent \textbf{Compare}  this paper recurrence (28) equivalent to (5) and proposition 2.2. in \cite[2009]{BSCS}  and note that (5) in \cite[2010]{BSCS} is just the same - as (58) in \cite[1878]{EdL} - the same except for notation.  The translation from "younger" notation of Bruce E. Sagan and Carla  D.  Savage (from one - left hand - side) into more matured by tradition  notation of Fran\c{c}ois  \'Edouard Anatole Lucas (from the other - right hand - side) is based on the identifications: the symbol used for $U$-binomials is  $\left\{... \right\}$  in place of  $\left(...\right)_U$  and 

$$ \left\{ n \right\} \equiv U_n \equiv n_{p,q}  , \ \left\langle n \right\rangle \equiv V_n.$$

\vspace{0.1cm}

\noindent For $s=t=1$  we get Fibonacci $U_n = F_n$ sequence with recurrence (41) becoming the recurrence known from Donald Ervin Knuth and Herbert Saul Wilf  masterpiece \cite[1989]{K-W}.  

\vspace{0.1cm}

\noindent \textbf{Example 3}. becomes \textbf{more general} case of the Theorem 2. if we allow $U$ to represent any zero characteristic field nonzero valued \textbf{functions'} sequence $U(x) = \left\langle U_n(x)\right\rangle_{n\geq 0}, \;  U_n(x) = \frac{p(x)^n-q(x)^n}{p(x)-q(x)} \equiv n_{p(x),q(x)},  p(x) \neq  q(x)$  i.e.  $p(x),q(x)$ denote distinct roots  of (20) and we have assumed  as well  that  $\frac{p(x)}{q(x)} $ are not roots of unity. 

\vspace{0.2cm}

\noindent \textbf{The End} of three examples.

\vspace{0.3cm}


\noindent Let us now come back to consider the case of $V$-binomials recurrence  according to what after Theorem 2a was called  -  "a not illusory obstacle".  
For that to do carefully let us at first make precise the main item. $ H(x) = \left\langle H_n(x)\right\rangle_{n\geq 0}$ is considered as a solution of \textcolor{red}{second} order recurrence (3) with peculiar case of (3) becoming the first order recurrence \textcolor{blue}{excluded}. This is equivalent to say that  $0\neq p \neq q \neq 0 $ and  $A \neq 0 \neq B$     (compare with \cite[1974]{Hilton 1974 partition}  by Anthony J. W. Hilton), where for notation convenience we shall again use awhile shortcuts  for  (3):  

$$a(x) \equiv a , b(x) \equiv b , s(x) \equiv s , t(x) \equiv t ,$$

$$p(x)\cdot q(x) = - t(x) \equiv p \cdot q = - t,$$

$$ H(x)= H(A(x),B(x),p(x),q(x)) = H(A,B) =A \cdot H(1,1;) ,$$
$$ H(1,1) = V(x)=V = (p-q)U(1,-1) , \  U = U(x)= U(1,1) , U(A,B) = U(A,B,p,q)$$
referring to this note (23) and the next to (23) abbreviations: 
$$H(x) \equiv H(a(x),b(x);s(x),t(x)) \equiv  H(A(x),B(x);s(x),t(x)),$$   $$H_n(x) \equiv H_n(a(x),b(x);s(x),t(x)),$$  
$$H_n(x) = H_n(A,B) = A\cdot p^n + B\cdot q^n \ H_n(1,1) = V_n(x)= V_n.$$

\vspace{0.2cm}

\noindent Now although $H(A,B) =(p-q)\cdot U(A,-B)$  the corresponding recurrences are different. Recall and then compare corresponding recurrences:

\noindent \textit{The identity }  
$$(p-q)\cdot(p^{r+s} - q^{r+s}) \equiv (p^{s+1} - q^{s+1})\cdot(p^r - q^r) - p\cdot q (p^{r-1} - q^{r-1}) \cdot(p^s - q^s)$$
due to $p\cdot q = -t$  is \textit{equivalent to }

$$
U_{r + s}(p,q) = U_{s+1} \cdot U_r(p,q) + t \cdot U_{r-1}\cdot U_s(p,q)
$$
(combinatorial derivation - see \cite[1999]{1BQ 1999},\cite[1999]{2BQ 1999}). This recurrence in its turn \noindent\textit{is equivalent to} 

$$
{r+s \choose r,s}_U = U_{s+1} \cdot {r+s-1 \choose r-1,s}_U + t \cdot U_{r-1}\cdot{r+s-1 \choose r,s-1}_U.
$$
\noindent Similarly \textit{the identity }  
$$(p-q)\cdot(p^{r+s} + q^{r+s}) \equiv (p^{s+1} - q^{s+1})\cdot(p^r + q^r) - p\cdot q (p^{r-1} + q^{r-1}) \cdot(p^s - q^s)$$
due to $p\cdot q = -t$  is \textit{equivalent to}

\begin{equation}
V_{r + s}(p,q) = U_{s+1} \cdot V_r(p,q) + t \cdot V_{r-1}\cdot U_s(p,q)
\end{equation}
For combinatorial interpretation derivation of the above for Fibonacci and Lucas sequences ($t=1=s$) see \cite[1999]{1BQ 1999}, \cite[1999]{2BQ 1999} and then see for more \cite[2003]{BQ2003} also by Arthur T. Benjamin and Jennifer J. Quinn.

\begin{defn}
Let  $ \left\{ n(x) \right\} \equiv U_n $ . $\left\langle n(x) \right\rangle \equiv V_n = p^n(x) +q^n(x)$, hence $V_0 = 2$ and  $V_1 = p + q = s(x)$. Let  $(p(x)- q(x)) \cdot U_n = p^n(x) - q^n(x)$, hence $U_0 = 0$ and  $V_1 = 1$; (roots are distinct). Then $V$-mixed-$U$  binomial coefficients for any $r,s \in \mathbb{N} \cup \{0\}$ ; are defined as follows
\begin{equation}
		{r+s \choose r,s}_{\left\langle . \right\rangle / \left\{ . \right\}} = \frac{V_{r+s}!}{V_r!\cdot U_s!}, 
\end{equation}
\end{defn}
${n \choose k}_{\left\langle . \right\rangle / \left\{ . \right\}}= 0 $  for $k>n$  and  ${n \choose 0}_{\left\langle . \right\rangle/\left\{ . \right\}}= 1.$

\vspace{0.2cm}

\noindent Note that:  ${r+s \choose r,s}_{\left\langle .\right\rangle/\left\{.\right\}} \neq {r+s \choose s,r}_{\left\langle .\right\rangle/\left\{.\right\}}.$

\vspace{0.2cm}

\noindent Note that: ${r+s \choose r,s}_{\left\langle .\right\rangle/\left\{.\right\}} \neq {r+s \choose r,s}_{ \left\{.\right\}/ \left\langle .\right\rangle }$

\vspace{0.3cm}

\noindent \textbf{\textcolor{red}{Theorem 2b}.} The recurrence (44) \noindent\textit{is equivalent to} 

\begin{equation}
{r+s \choose r,s}_{\left\langle . \right\rangle / \left\{ . \right\}} = U_{s+1} \cdot {r+s -1  \choose r-1,s}_{\left\langle . \right\rangle / \left\{ . \right\}} + t \cdot U_s \cdot {r+s-1 \choose r,s-1}_{\left\langle . \right\rangle / \left\{ . \right\}}.
\end{equation}

\vspace{0.5cm}

\noindent \textbf{Example 4.} Application of  Theorems 1, 2a and  the Theorem 2b method.  
\noindent Recall the characteristic equation notation:  $z^2 =\textcolor{blue}{z^2 = s \cdot x + t}$ as is common in many publications with the restriction:  $p,q$ roots are distinct. Recall that $p+q=s$  and $p\cdot q=-t$. Let $\Delta = s^2 - 4t$. Then $\Delta = (p-q)^2$. Hence we have as in \cite[1878]{EdL}

$$
2\cdot U_{r+s} =  U_r V_s  +  U_s V_r ,\ \ \ \ 
2\cdot V_{r+s} =  V_r V_s  +  \Delta \cdot U_s U_r. 
$$
Consider   $\textbf{V}$-\textcolor{blue}{addition formula} i.e. the second  of two trigonometric-like addition formulas (42) from \cite[1878]{EdL}. One readily recognizes that  \textit{according  to the Theorem 2b  the following is true}:     $2 \cdot V_{r+s} =  V_r V_s  +  \Delta \cdot U_s U_r$  \textit{is equivalent to}

\begin{equation}
2\cdot {r+s \choose r,s}_{\left\langle . \right\rangle / \left\{ . \right\}} = V_s \cdot {r+s -1  \choose r-1,s}_{\left\langle . \right\rangle / \left\{ . \right\}} + \Delta \cdot U_r \cdot {r+s-1 \choose r,s-1}_{\left\langle . \right\rangle / \left\{ . \right\}}.
\end{equation}
The combinatorial interpretation of mixed binomials is  part of the subject of a forthcoming "in statu nascendi"  note.

\vspace{0.6cm}

\section{Snatchy information on combinatorial interpretations of $H(x=1)$-binomials and their relatives.}

\vspace{0.1cm}

\noindent \textbf{4.1.} 
\noindent In regard to  \textbf{combinatorial interpretations} of $H$-\textbf{bi}nomial or $F$-\textbf{multi}nomial coefficients or related arrays we leave that subject apart from this note.  Nevertheless we direct the reader  to some comprise papers and references therein  via listing; these are here \textsl{for example}  the following:

\vspace{0.2cm}

\noindent Listing. \textbf{1.} \cite[1984]{Voigt} by Bernd Voigt: on common generalization of binomial coefficients, Stirling numbers and Gaussian coefficients .

\vspace{0.2cm}

\noindent Listing. \textbf{2.}  \cite[1991]{wachs} by Michelle L. Wachs and   Dennis White and in \cite[1994]{wachs 2} by Michelle L. Wachs: on p,q-Stirling numbers and set partitions.

\vspace{0.2cm}

\noindent Listing. \textbf{3.} \cite[1993]{medicis} by Anne De M\'edicis  and Pierre Leroux:  on Generalized Stirling Numbers, Convolution Formulae and (p,q)-Analogues.

\vspace{0.2cm}

\noindent Listing. \textbf{4.} \cite[1998]{Konva 1998} John Konvalina:  on generalized binomial coefficients and the  Subset-Subspace Problem.  Consult examples [Ex. q* ; 6] and  [Ex. q* ; 7] in \textbf{4.3}. below.  Then see also the article \cite[2000]{Konva 2000}  by   John Konvalina on an unified  simultaneous  interpretation of binomial coefficients of both kinds, Stirling numbers of both kinds and Gaussian binomial coefficients  of both kinds.

\vspace{0.2cm}

\noindent Listing. \textbf{5.} Ira  M. Gessel  and Xavier G\'erard Viennot in  \cite[1989]{G-V} deliver now the well known their interpretation of the fibonomials in terms of non-intersecting lattice 
paths.

\vspace{0.2cm}

\noindent Listing. \textbf{6.} In \cite[2004]{RW}  Jeffrey Brian Remmel and Michelle L. Wachs derive a new rook theory interpretation of a certain class of generalized
Stirling numbers and their $(p,q)$-analogues.  In particular they prove that their $(p,q)$-analogues of the generalized Stirling numbers of the second
kind may be interpreted in terms of colored set partitions and colored restricted growth
functions.

\vspace{0.2cm}

\noindent Listing. In \textbf{7.}  \cite[2005]{Otta vio Munari}  by Ottavio M. D'Antona  and Emanuele Munarini deal - in terms of weighted binary paths -  with   combinatorial interpretation of the connection constants which is in particular unified, simultaneous combinatorial interpretation  for Gaussian coefficients, Lagrange sum, Lah numbers, ,  q-Lah numbers, Stirling numbers of both kinds , q-Stirling numbers of both kinds. Note the usefull correspondence: weighted binary paths $\Leftrightarrow$     edge  colored binary paths

\vspace{0.2cm}
\noindent  Listing. \textbf{8.} Maciej  Dziemia\'nczuk in \cite[2011]{MD 2011}  extends the results of John Konvalina from \textbf{4.} above. The Dziemia\'nczuk'  $\zeta$ - analogues of the Stirling numbers arrays of both kinds cover ordinary binomial and Gaussian coefficients, $p,q$-Stirling numbers and other combinatorial numbers studied with the help of object selection, Ferrers diagrams and rook theory. The $p,q$-\textbf{binomial} arrays are special cases of $\zeta$- numbers' arrays, too. 

\vspace{0.1cm}
\noindent $\zeta$ -number of the  first and the second kind is the number of ways to select $k$  objects from $k$  of $ n$ boxes without box repetition allowed and with box repetition allowed, respectively.  The weight vectors used for objects constructions and  statements derivation are functions of parameter $\zeta$.


\vspace{0.2cm}
\noindent  Listing. \textbf{9.} As regards combinatorial interpretations via  tilings in   \cite[2003]{BQ2003} and  \cite[2010]{BP}  - see \textbf{4.2.} below.  


\vspace{0.2cm}

\noindent Listing. \textbf{10.} In   \cite[2003]{Cigler} Johann Cigler  introduces "abstract Fibonacci polynomials"  - interpreted in terms of Morse coding sequences monoid with concatenation (monominos and dominos tiling then).  Cigler's abstract Fibonacci polynomials are monoid algebra over reals valued polynomials with straightforward Morse sequences i.e. tiling recurrence originated (1.6)  "addition formula"

$$F_{m+n}(a,b) =  F_{m+1}(a,b)\cdot F_{m}(a,b) + b\cdot F_{n-1}(a,b)\cdot F_{n}(a,b),$$
which is attractive and seductive to deal with   within the context of this paper Theorem 1. The combinatorial tiling interpretation of the model is its construction framed in the Morse coding sequences monoid with concatenation (monominos and dominos tiling then).


\vspace{0.2cm}

\noindent  Listing. \textbf{11.}...to be juxtaposed to Listing. \textbf{12.}.

\noindent In \cite[2009]{Cigler q-Lucas 2010} Johann Cigler considers special Ward-Horadam   \textit{polynomials  sequences} and reveals the tiling combinatorial interpretation of these special Ward-Horadam   \textit{polynomials  sequences} in the spirit of Morse with monomino, domino alphabet monoid as here above in Listing. \textbf{10.}. Compare with technique in
\cite[1994]{Foata 1994 Pisa}  by Dominique Foata  and  Guo Niu Han  \textbf{where }one considers "de rubans d'ordre  $m$, qu'on peut repr\'esenter comme $m$  rubans remplis de monominos et de dominos, de m\'eme longueur $n$"  and where  Pafnuty Lvovich Chebyshev polynomials are treated.

\noindent \textbf{1.} In \cite[2009]{Cigler q-Lucas 2010}   by  Johann Cigler the $q$-Fibonacci polynomial $F_n(x,s,q) = \sum_{c\in \Phi_n}w(c)\equiv w(\Phi_n)$ is the $q$-weight function of the set $\Phi_n$ of all words  (coverings) $c$ of length $n-1$ in Morse (tiling) alphabet $\left\{a,b\right\}$ i.e.it is the  corresponding generation function for number of linear \textcolor{red}{\textbf{$q$-weighted}} tilings as clearly  $\Phi_n$  may be identified with the set of all linear \textcolor{red}{\textbf{$q$-weighted}} tilings  of  $(n-1)\times 1$  rectangle or equivalently with Morse code sequences of length $n-1$. 

\noindent Then classical Fibonacci polynomials $F_n(x,s,1)$ satisfy  this paper recursion (3) i.e the Ward-Horadam functions' sequence recursion below -  with $H_0(x)=0$ , $H_1(x)=x$ ; $s(x)=x$ and $t(x)=s$; ( here $H \mapsto F$):

$$H_{n+2}(x) = s(x) \cdot H_{n+1}(x) + t(x)\cdot H_n, \; n\geq 0, \ H_0 = a(x), H_1 = b(x).$$

\noindent The $F(x,s,q)$-binomial array  ${n \choose k}_{F(x,s,q)} $  is not considered in  \cite[2009]{Cigler q-Lucas 2010}. Similarly:
\vspace{0.1cm}

\noindent \textbf{2.} the $q$-Lucas polynomial $L_n(x,s,q) = \sum_{c\in \Lambda_n}w(c) \equiv w(\Lambda)$ is the $q$-weight function of the set $\Lambda_n$ of  all coverings $c$ with arc monominos and dominos of the circle whose circumference has length $n$. Hence $L_n(x,s,q)$ is corresponding generation function for number of \textcolor{blue}{\textbf{$q$-weighted}} tilings of the circle whose circumference has length $n$. It may be then  combinatorially seen that $w(\Lambda_n)= w(\Phi_{n+1}) + s\cdot w(\Phi_{n-1}))$  hence specifically for the classical Fibonacci and Lucas polynomials we get the well known relation $L_n(x,s,1) = F_{n+1}(x,s,1)+ s\cdot F_{n-1}(x,s,1)$.  In  \cite[2009]{Cigler q-Lucas 2010} Johann Cigler proves more : $L_n(x,s,q) = F_{n+1}(x,s,q)+ s\cdot F_{n-1}(x,s,q)$ and provides a combinatorial interpretation of this relation, too.

\noindent Polynomials $L_n(x,s,1)$ satisfy this paper recursion (3)  with $H_0(x)=2$ , $H_1(x)=x$ ; $s(x)=x$ and $t(x)=s$; (here $H \mapsto L$)

\noindent The $L(x,s,q)$-binomial array  ${n \choose k}_{L(x,s,q)} $  is not considered in  \cite[2009]{Cigler q-Lucas 2010}.


\vspace{0.2cm}

\noindent Listing. \textbf{12.} ..to be juxtaposed to Listing. \textbf{11.}.

\noindent \cite[2009]{BSCS}   by  Bruce E. Sagan  and Carla  D.  Savage  the symbol $\left\{n\right\}\equiv U_n$    denotes the $n-th$ element of the  fundamental Lucas sequence $U$ satisfying this paper recurrence (2) with initial conditions $\left\{0\right\}=0$, $\left\{1\right\}=1$.  Naturally  $\left\{n\right\}$ is a polynomial in parameters  $s,t$. So is also the $U$-binomial coefficient $\left\{ n \atop k\right\}_U \equiv \left\{ n \atop k\right\}_{p,q}$. 

\vspace{0.1cm}

\noindent Similarly -  the symbol $\left\langle n \right\rangle\equiv V_n$    denotes the $n-th$ element of the  primordial Lucas sequence $V$ satisfying this paper recurrence (2) with initial conditions $\left\langle 0 \right\rangle = 2$ , $\left\langle 1 \right\rangle = s$ .  Naturally  $\left\langle n \right\rangle $  is a polynomial in parameters  $s,t$. So is also the $V$-binomial coefficient $\left\{ n \atop k\right\}_V \equiv \left\langle  n \atop k\right\rangle_{p,q}$.  $V$-binomials are not considered in \cite[2009]{BSCS}. Both fundamental and primordial sequences are interpreted via tilings similarly to the above in  \textbf{11.}  Johann Cigler attitude  rooted in already text-books tradition - see for example \cite[1989]{GKP 1989}  by Ronald Graham, Donald Ervin  Knuth  and Oren Patashnik. 

\vspace{0.1cm}

\noindent An so:  $\left\{n\right\}$ is generation function for number of linear tilings  of  $(n-1)\times 1$  rectangle or equivalently of number of Morse code sequences of length $n-1$.

\vspace{0.1cm}

\noindent $\left\langle n \right\rangle $ is generation function for number of circular tilings of the circle whose circumference has length $n$. Using naturally proved (just seen) relations  Bruce E. Sagan  and Carla  D.  Savage  derive \textbf{two} combinatorial interpretations of the the same $\left\{ m+n \atop m,n\right\}_{p,q}$  via  Theorem 3.1. from which we infer the following.

\vspace{0.1cm}

\noindent \textbf{ 1.}   $\left\{ m+n \atop m,n\right\}_{p,q}$ is  the weight of all linear tilings of all integer partitions  $\lambda$ inside the $m \cdot n$ rectangle 

\noindent hence $\left\{ m+n \atop m,n\right\}_{p,q}$ is the generating function  for numbers of such tilings of partitions.

\vspace{0.1cm}

\noindent  \textbf{2.}  $ 2^{m+n}\cdot \left\{ m+n \atop m,n\right\}_{p,q}$  is  the weight of all circular tilings of all integer partitions  $\lambda$ inside the $m \cdot n$ rectangle 

\noindent hence $\left\{ m+n \atop m,n\right\}_{p,q}$ is the generating function  for numbers of such tilings of partitions.

\vspace{0.1cm}

\noindent \textbf{Explanation.} from \cite[2009]{BSCS}. \textit{A linear tiling of a partition}  $\lambda$   is a covering of its Ferrers diagram with disjoint dominos and monominos obtained by linearly tiling each  $\lambda_i$ part.  In circular tiling of a partition  $\lambda$ one performs circular tiling of each  $\lambda_i$ part


\vspace{0.2cm}

\noindent Listing. \textbf{13.} In  \cite[2009]{MD Counting Bipartite 2009} \textit{Counting Bipartite, $k$-Colored and Directed Acyclic Multi Graphs Through $F$-binomial coefficients by Maciej Dziemia\'nczuk}, the explicit relation between the number of the so called in  \cite[2009]{MD Counting Bipartite 2009} $k$-colored $\alpha$-multigraphs and  $N(\alpha)$-multinomial coefficients has been established.   

\noindent Here it is a directed acyclic graph with $\alpha$-multiple edges (i.e. any two vertices might be connected by at most $(\alpha-1)$ directed edges ) what is called acyclic $\alpha$-multi digraph.

\vspace{0.1cm}

\noindent Let $N(\alpha) =  \{n_{\alpha,\alpha}\}_{n\geq 0}$   i.e.  one considers the  $p=q \equiv \alpha$ case of $U = \left\langle n_{p,q}\right\rangle$ number sequence, hence 
$n_{\alpha,\alpha} =  n \cdot \alpha^{n-1}$. For example (we quote these  from \cite[2009]{MD Counting Bipartite 2009}):

\begin{enumerate}
\item $N(1) = 0, 1, 2, 3, 4, 5, 6, 7, 8, 9, 10, ... $
\item $N(2) = 0, 1, 4, 12, 32, 80, 192, 448, 1024, 2304, 5120, ...  $
\item $N(3) = 0, 1, 6, 27, 108, 405, 1458, 5103, 17496, 59049, 196830, ...  $
\item $N(4) = 0, 1, 8, 48, 256, 1280, 6144, 28672, 131072, 589824, 2621440,...   $
\end{enumerate}

\vspace{0.1cm}

\noindent Maciej Dziemia\'nczuk proves in  \cite[2009]{MD Counting Bipartite 2009} that \textbf{if} $\gamma_{\alpha,n,k}$ is the number of all $k$-colored $\alpha$-multigraphs with $n$ vertices, \textbf{then} 

$$
	\gamma_{\alpha,n,k} = \sum_{{b_1+...+b_k=n \atop b_1,...,b_k\geq 0}} \fnomialF{n}{b_1,b_2,...,b_k}{N(\alpha)}
$$
The case of $k$-colored graphs \textit{without multiple edges} i.e. $\alpha=2$ -  is to be found in \cite[2003]{finch}  by  Steven R. Finch.

\vspace{0.1cm}
\noindent Maciej Dziemia\'nczuk in  \cite[2009]{MD Counting Bipartite 2009} proves also that the unsigned values of the first row of inversion matrix for $N(\alpha)$-binomial coefficients considered there are equal to the numbers of directed acyclic $\alpha$-multigraphs with $n$ nodes. 

\vspace{0.1cm}

\noindent \textbf{Ad Notation Remark.}  The author of \cite[2009]{MD Counting Bipartite 2009} and  \cite [2008]{MD Inv F-binomial 2008} - uses misleading, cumbersome notation. He writes $m_{n,k}\:^{-1} =\fnomial{n}{k}^{-1}$ instead of $ M^{-1}\:_{n.k} =  \left( \fnomial{r}{s} \right)^{-1}\:_{n,k}$ , where here $M= \left(m_{n,k} \right) $.

\vspace{0.1cm}

\noindent The exemplification that now  follows is concerned with the referring to Steven R. Finch particular result from \cite[2009]{MD Counting Bipartite 2009}. For that it might be helpful  to  consult eventually also \textbf{$\psi$-mulinomial remark} in this note Section 2.

\noindent In view of the final Remark in \cite[2004]{akk ASCM 2004},  the fundamental logarithmic Fib-Binomial Formula (\cite{akk ASCM 2004}, Section 4)
$$
	\phi_n^{(t)}(x +_F a) \equiv \left[ \mathrm{exp}\{a\partial_F\}\phi_n^{(t)} \right](x)
	= \sum_{k\geq 0}\left[\!\! \begin{array}{c} n \\ k \end{array}\!\!\right]_F \phi_{n-k}^{(t)}(a)x^k
$$
$$
	t=0,1; \ \ \ |x| < a; \ \ \ n\in \mathbb{Z}
$$
\noindent may be considered as $F$-Binomial formula for any natural numbers valued sequence $F$ with $ F_0 =1$  (the class considered in \cite[2004]{akk ASCM 2004} is much broader). 

\vspace{0.2cm}

\noindent Now put in $F$-binomial formula above $t=0$ and $a=x=1$ and pay attention to - that according to the Steven M. Roman definition of hybrid binomial coefficients in 
\cite[1992]{Steven M. Roman 1992} one states -  this time for $F$-hybrid binomial coefficients - the wanted for  $n,k\geq 0$ identification  
$$
	\left[\!\! \begin{array}{c} n \\ k \end{array}\!\!\right]_F = \fnomial{n}{k}.
$$

\vspace{0.2cm}

\noindent For special $F$-sequences introduced  as $F$-cobweb posets admissible sequences  (\cite[2003]{akk ITSF 2003}, \cite[2001]{ akk ITSF 2001}, \cite [2005]{akk Allahabad 2005}) the $F$-binomial coefficients for $n,k\geq 0$ acquire joint combinatorial interpretation from Listing. \textbf{14.} 

\noindent Then in Morgan Ward  Calculus of Sequences (\cite[1936]{Ward}) notation  adapted in umbra calculus by Kwa\'sniewski  one may write the following appealing formula :
$$
	\left( 1 +_F 1 \right)^n \equiv \sum_{n\geq 0}\fnomial{n}{k}
$$
which is a special case of 
$$
	\left( 1 +_F... +_F 1 \right)^n \equiv \sum_{n\geq 0}\fnomial{n}{k}.
$$
In particular application  ( see  \cite [2009]{MD Counting Bipartite 2009} by  Maciej  Dziemia\'nczuk) we may then write in Ward Morgan \cite[1936]{Ward} spirit notation that the number $\gamma_{n,2}$ of all 2-colored graphs as in \cite{finch} is now due to \cite [2009]{MD Counting Bipartite 2009} equal to  
$$
	\gamma_{n,2} = \sum_{k\geq 0} {n\choose k} \cdot 2^{k(n-k)} = \sum_{k\geq 0} \fnomialF{n}{k}{N(2)} \equiv \left( 1 +_{N(2)} 1 \right)^n
$$
For more - consult  \cite [2009]{MD Counting Bipartite 2009} by  Maciej Dziemia\'nczuk.


\vspace{0.5cm}

\noindent Listing. \textbf{14.} \textbf{Cobweb posets' partitions and hyper-boxes tilings.} 

\vspace{0.2cm}

\noindent \textbf{14.1}  \textit{A definition of the cobweb poset}.

\vspace{0.2cm}

\noindent The cobweb posets where introduced under this name  in several paper - see: \cite[2007]{11} , \cite{2}-\cite{6} ,  \cite[2008]{10}  and references therein.  There - the  cobweb posets where defined  in terms of their poset Hasse diagrams.   The definition we are here delivering is taken from  \cite[2009]{91}, \cite[2009]{2}  (for equivalent definitions see also: \cite[2009]{8}, \cite[2009]{9},  \cite[2009]{7}),  \cite[2010]{91}).   

\noindent Namely; let us consider any infinite  chain  $\left\{\Phi_k\right\}_{k\geq 0}$  of  trivial posets  $\equiv$ antichains  $\equiv$ a chain of  independent sets  $\equiv$ a  chain of  trivial unary relations, then we define

\begin{defn} (cobweb poset $\Pi_n$) 
\noindent Let  $n\in N \cup \left\{0\right\}\cup \left\{\infty\right\}$. Then
$$\Pi_n =  \oplus_{s= 0}^n \Phi_s ,$$
$$\Pi =  \oplus_{s\geq 0} \Phi_s, $$
where $\oplus$ denotes ordinal sum of posets. 
\end{defn}
\vspace{0.1cm}
\noindent  Recall that the ordinal sum [linear sum] of two disjoint ordered sets $P$ and $Q$, denoted by  $P\oplus Q$, is the union of $P$ and $Q$, with $P$'s elements ordered as in $P$ while $Q$'s elements are correspondingly ordered as in $Q$, and  for each $x\in P$ and $y \in Q$ we put  $x \leq y$ .  
\noindent  The Hasse diagram of  $P\oplus Q$ we construct placing $Q$'s diagram just above $P$'s diagram and with an edge between each minimal element of $Q$ and each maximal element of $P$. 

\vspace{0.1cm}
\noindent  The cobweb posets might be identified with a chain of di-bicliques i.e. by definition - a chain of complete bipartite one direction digraphs; see for example \cite[2009]{91}, \cite[2009]{2}, \cite[2009]{8}, \cite[2008]{10}. Any chain of relations is therefore obtainable from the cobweb poset chain of complete relations via deleting arcs in di-bicliques of the complete relations chain.

\begin{figure}[ht]
\begin{center}
	\includegraphics[width=75mm]{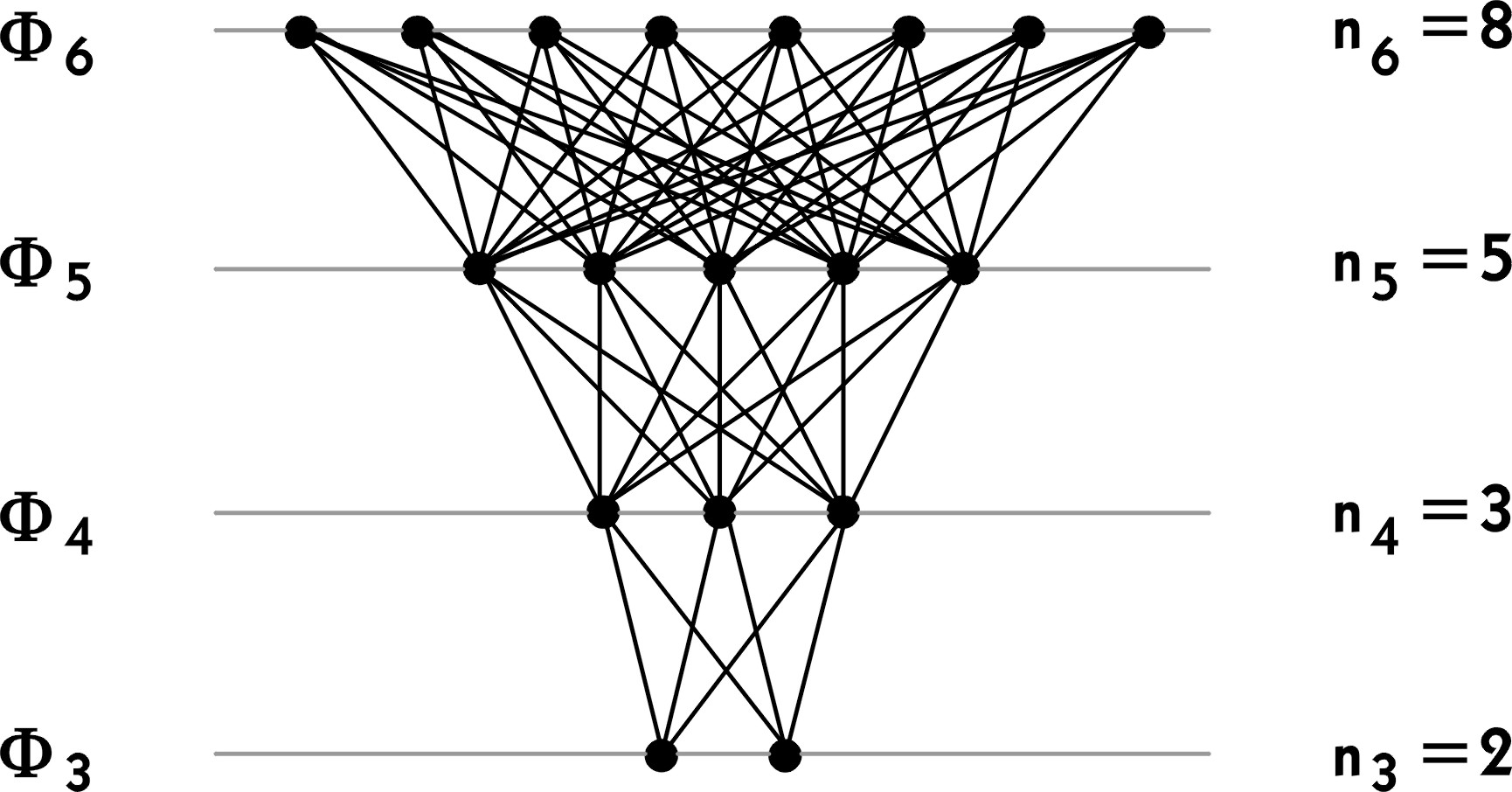}
	\caption{Display of four levels of Fibonacci numbers' finite Cobweb sub-poset }
\end{center}
\end{figure}

\begin{figure}[ht]
\begin{center}
	\includegraphics[width=75mm]{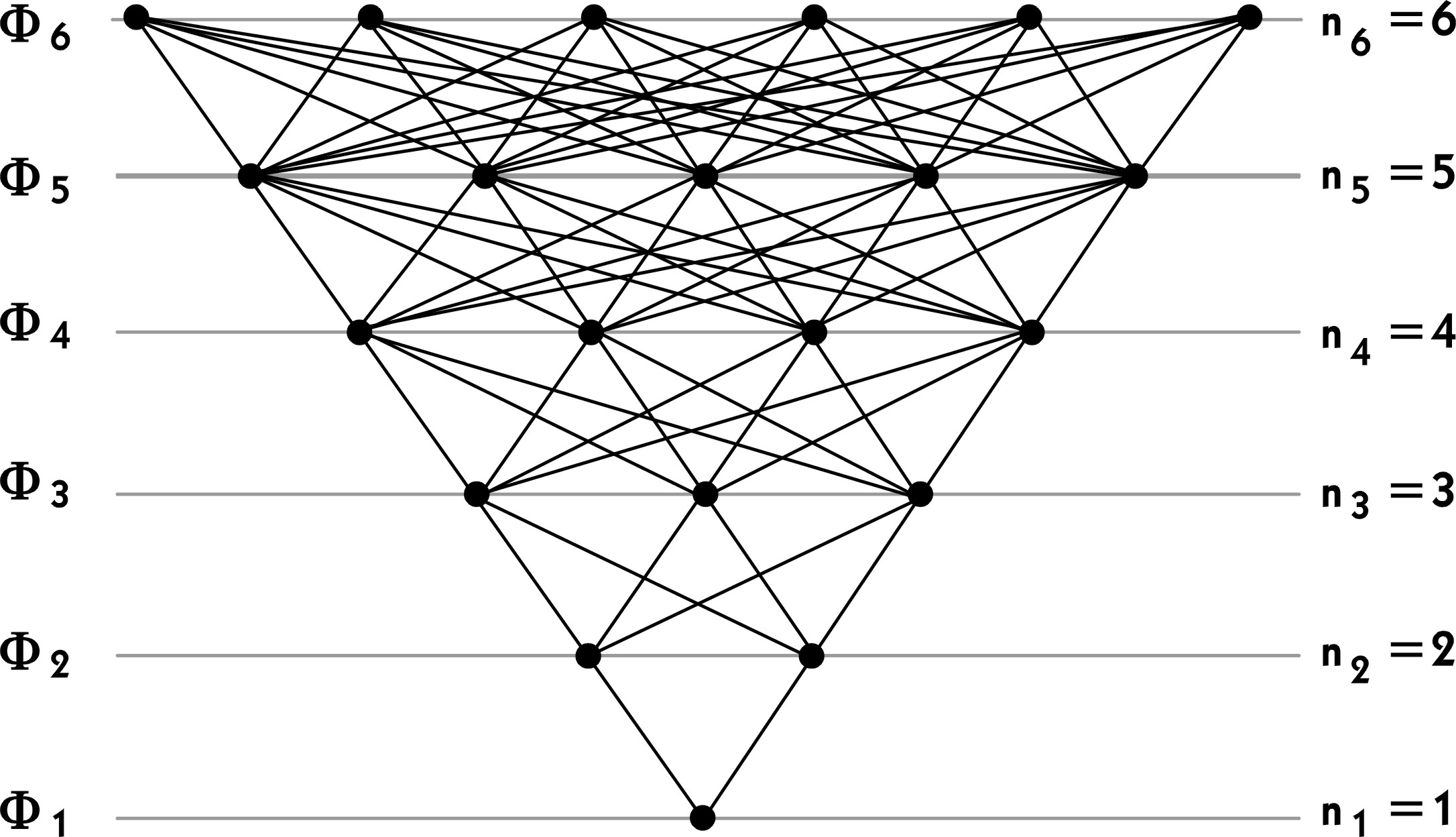}
	\caption{Display of Natural numbers' finite Cobweb sub-poset }
\end{center}
\end{figure}



\vspace{0.2cm}

\noindent \textbf{14.2}  \textbf{Combinatorial interpretation of $F$-binomials in terms of cobweb poset}.

\vspace{0.1cm}

\noindent Combinatorial interpretation of cobweb posets via their cover relation digraphs (Hasse diagrams) called KoDAGs has been developed in  \cite[2009]{8},  \cite[2010]{6}, \cite[2009]{91}, \cite[2009]{2}, \cite[2010]{3}, \cite[2010]{4}. The recent equivalent formulation of this combinatorial interpretation is to be found also in \cite[2009]{7} or \cite[2009]{9} from which we quote it here down. We shall now start to use the upside down notation :  $F_n \equiv n_F$  as an extension of $n_{p,q} \equiv U_n(p,q)$  or specifically $n_{1,q} \equiv U_n(1,q) \equiv n_q$ notation.

\begin{defn}
Admissible sequence $F$-\textbf{binomial} \textbf{coefficients} are defined as follows
$$
	\fnomial{n}{k} = \frac{n_F!}{k_F!(n-k)_F!} 
	= \frac{n_F\cdot(n-1)_F\cdot ...\cdot(n-k+1)_F}{1_F\cdot 2_F\cdot ... \cdot k_F}
	= \frac{n^{\underline{k}}_F}{k_F!}
$$
\noindent while $n,k\in \mathbb{N}$ and $0_F! = n^{\underline{0}}_F = 1$  with $n^{\underline{k}}_F \equiv \frac{n_F!}{k_F!}$ staying for falling factorial.
Zero or natural numbers' valued $F$ sequence is called  $F$-graded poset  \textcolor{blue}{\textbf{admissible}} sequence iff  $\fnomial{n}{k} \in N \cup\left\{0\right\}$ ( In particular we shall use the expression - $F$-cobweb admissible sequence).
\end{defn}

\begin{defn}
$$
C_{max}(\Pi_n) \equiv  \left\{c=<x_0,x_1,...,x_n>, \: x_s \in \Phi_s, \:s=0,...,n \right\} 
$$  
i.e. $C_{max}(\Pi_n)$ is the set of all maximal chains of $\Pi_n$
\end{defn}

\noindent and consequently (see Section 2 in \cite[2011]{12}  on Cobweb posets' coding via $N^\infty$ lattice boxes)

\begin{defn} ($C^{k,n}_{max}$)  

$$
C_{max}\langle\Phi_k \to \Phi_n \rangle \equiv \left\{c=<x_k,x_{k+1},...,x_n>, \: x_s \in \Phi_s, \:s=k,...,n \right\}\equiv
$$

$$
	\equiv \big\{ \mathrm{maximal\ chains\ in\ } \langle \Phi_k \rightarrow \Phi_n \rangle \big\} \equiv
	C_{max}\big( \langle \Phi_k \rightarrow \Phi_n \rangle \big) \equiv
	C^{k,n}_{max}.
$$
\end{defn}

\noindent \textbf{Note.} The $C_{max}\langle\Phi_k \to \Phi_n \rangle \equiv C^{k,n}_{max}$
is the hyper-box points'  set \cite{12}, \cite{2} of  Hasse sub-diagram corresponding maximal chains and it defines biunivoquely 
the layer $\langle\Phi_k \to \Phi_n \rangle = \bigcup_{s=k}^n\Phi_s$  as the set of maximal chains' nodes (and vice versa) -
for  these arbitrary $F$-denominated \textbf{graded} DAGs (KoDAGs included).

\vspace{0.1cm}

\noindent The formulation of the fractals reminiscent combinatorial interpretation of cobweb posets via their cover relation digraphs (Hasse diagrams) is the following observation  (named deservedly ? a Theorem) equivalent to that of \cite{8}, \cite{7},  \cite{9}.

\vspace{0.2cm}

\noindent \textbf{Theorem 3} \cite{11}, \cite{9}, \cite{8}, \cite{3} \\
\noindent \textit{For $F$-cobweb admissible sequences $F$-binomial coefficient $\fnomial{n}{k}$ is the cardinality of the family of \emph{equipotent} to  $C_{max}(P_m)$ mutually disjoint maximal chains sets, all together \textbf{partitioning } the set of maximal chains  $C_{max}\langle\Phi_{k+1} \to \Phi_n \rangle$  of the layer   $\langle\Phi_{k+1} \to \Phi_n \rangle$, where $m=n-k$.}

\vspace{0.1cm}

\noindent \textbf{Comment 1}. For the  above  combinatorial  interpretation of  $F$-binomials' array \textit{it does not matter}  of course whether the diagram is being directed  or not, as this combinatorial interpretation is  equally valid for partitions  of the family of  $SimplePath_{max}(\Phi_k - \Phi_n)$ in  comparability graph of the Hasse  digraph with self-explanatory notation used on the way. The other insight into this irrelevance for combinatoric interpretation is in \cite{9}: colligate the coding of $C^{k,n}_{max}$ by hyper-boxes.

\begin{figure}[ht]
\begin{center}
		\includegraphics[width=75mm]{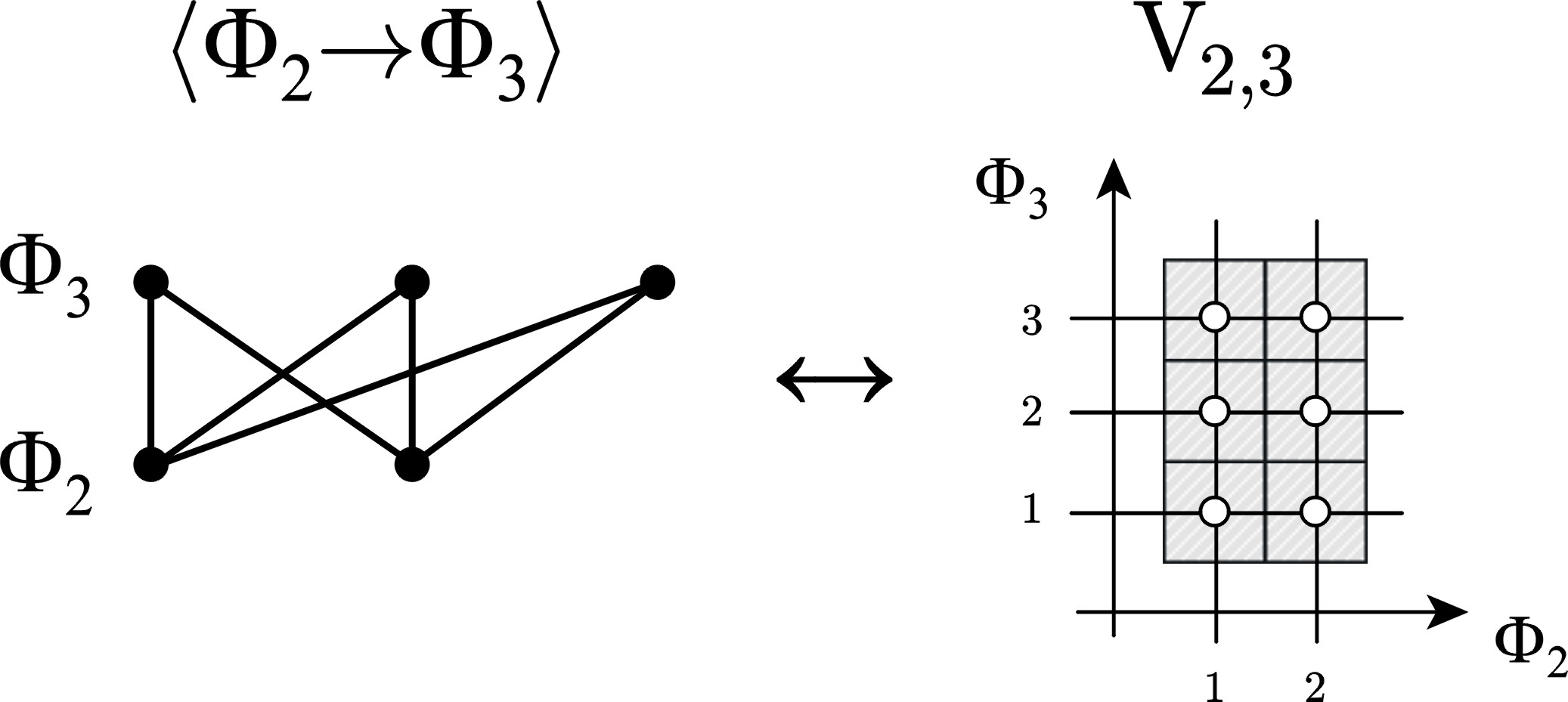}
	\caption{Bipartite  layer $\langle\Phi_3 \rightarrow \Phi_4 \rangle$ with six maximal chains and equivalent hyper-box $V_{2,3}$  with six white circle-dots} \label{fig:representation}
\end{center}
\end{figure}

\begin{figure}[ht]
\begin{center}
	\includegraphics[width=75mm]{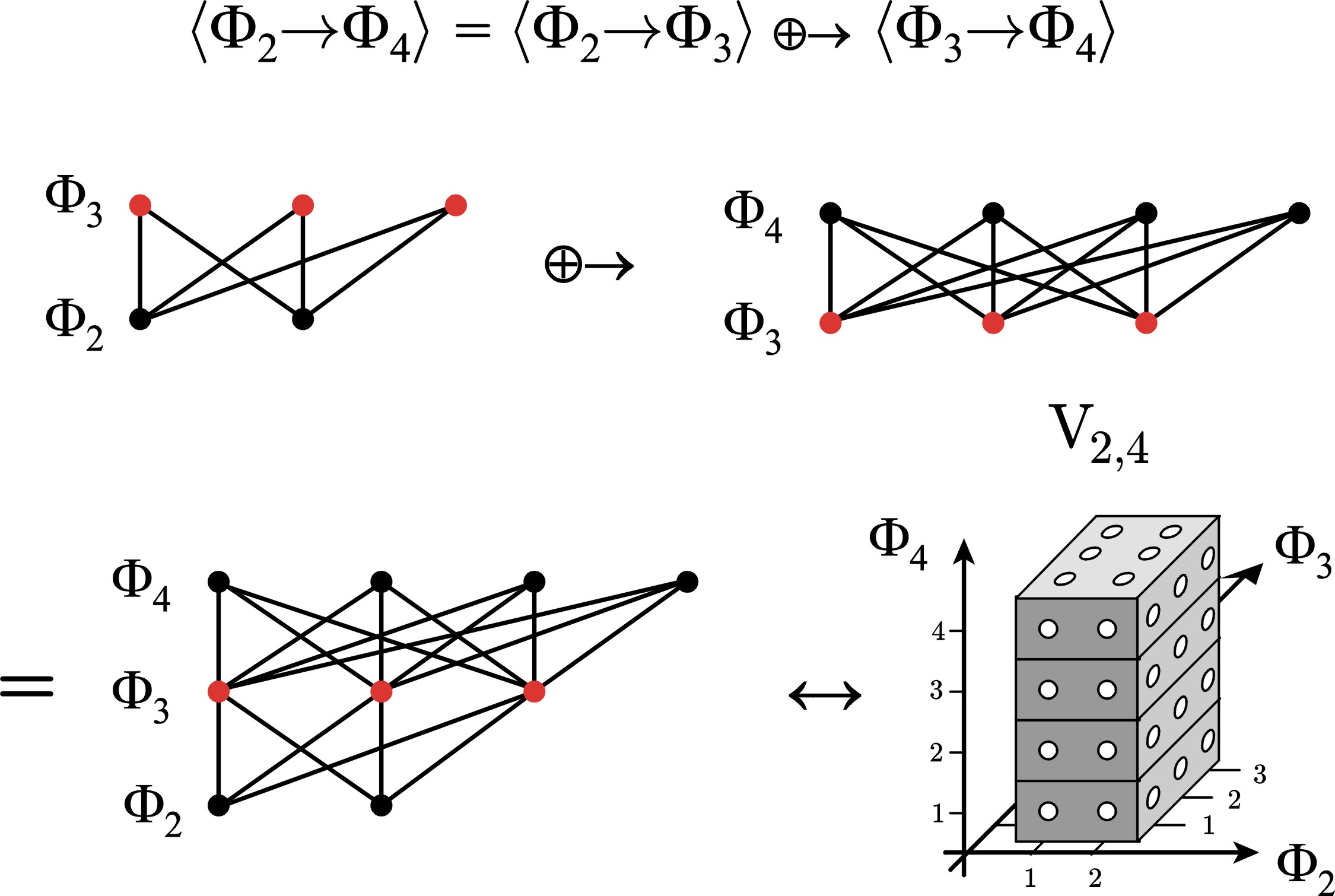}
	\caption 	{Display of the natural join of bipartite layers $\left\langle \Phi_k \rightarrow \Phi_{k+1} \right\rangle$ $F = N $, resulting in  $2\cdot3\cdot4$ maximal chains and equivalent hyper-box $V_{2,4}$ with $2\cdot3\cdot4$ white circle-dots.}\label{fig:representation} 
\end{center}
\end{figure}

\vspace{0.2cm}

\noindent \textbf{14.3}  \textbf{Combinatorial interpretation of $F$-binomials in discrete hyper-boxes language}.

\vspace{0.1cm}

\noindent  In order to formulate the above combinatorial interpretation let us recall indispensable identifications   from  \cite[2009]{3},  \cite[2009]{9}. 

\vspace{0.2cm}

\noindent $C_{max}(\Pi_n)$ is the set of all maximal chains of $\Pi_n$.
 
$$C^{k,n}_{max} = \big\{ \mathrm{maximal\ chains\ in\ } \langle \Phi_k \rightarrow \Phi_n \rangle \big\}.$$

\noindent Consult  Section 3. in  \cite[2011]{12} in order to view $C_{max}(\Pi_n)$ or  $C^{k,n}_{max}$ as the hyper-boxes of points.\\ 
Now following \cite[2009]{3}, \cite[2009]{2}), \cite[2009]{91}, \cite[2011]{12}) - we define the $F$-hyper-box $V_{k,n}$.

\begin{defn}  The  the discrete finite rectangular $F$-hyper-box or $(k,n)-F$-hyper-box or in everyday parlance just $(k,n)$-box  $V_{k,n}$ is the Cartesian product
$$
	V_{k,n} = [k_F]\times [(k+1)_F]\times ... \times[n_F].
$$
\end{defn}

\noindent We identify the following two object just by agreement according to the $F$-natural identification:
$$
	C^{k,n}_{max} \equiv V_{k,n}
$$
i.e.
$$
C^{k,n}_{max} = \big\{ \mathrm{maximal\ chains\ in\ } \langle \Phi_k \rightarrow \Phi_n \rangle \big\} \equiv V_{k,n}.
$$
\vspace{0.1cm}

\noindent For numerous illustrations via natural join of posets alike that by Fig.4. -  see \cite[2009]{2}), \cite[2009]{91}  and \cite[2011]{12}.

\vspace{0.1cm}
 
\noindent Now in  discrete hyper-boxes language the \textit{combinatorial interpretation reads}:

\vspace{0.1cm}

\noindent \textbf{Theorem 4}. \cite[2009]{91} \\
\noindent \textit{For $F$-cobweb admissible sequences $F$-binomial coefficient} $\fnomial{n}{k}$ \textit{is the cardinality of the family of \emph{equipotent} to}  $V_{0,m}$ \textit{mutually disjoint 
discrete hyper-boxes, all together \textbf{partitioning } the discrete hyper-box }  $V_{k+1,n}$    $\equiv$   \textit{the layer}   $\langle\Phi_{k+1} \to \Phi_n \rangle$,\textit{ where} $m=n-k$.

\vspace{0.3cm}

\noindent \textbf{14.4}  \textbf{\textcolor{blue}{The cobweb tiling problem} in the language of \textcolor{blue}{discrete hyper-boxes}}.

\vspace{0.2cm}

\noindent \textbf{Comment 2.} General  "fractal-reminiscent" comment. The discrete $m$-dimensional $F$-box ($m = n-k$) with edges' sizes designated by natural numbers' valued sequence $F$ 
where invented in \cite[2010]{12} as a response to the so called \emph{ cobweb tiling problem} posed in \cite[2007]{11} and then repeated in \cite[2009]{8}. This tiling problem was considered by Maciej Dziemia\'nczuk in \cite[2008]{13} where it was shown that not all admissible $F$-sequences permit tiling as defined in \cite[2007]{11}.  Then - after \cite[2009]{12} this tiling  problem was considered by  Maciej Dziemia\'nczuk in discrete hyper-boxes language \cite[2009]{48}.

\vspace{0.1cm}

\begin{figure}[ht]
\begin{center}
	\includegraphics[width=75mm]{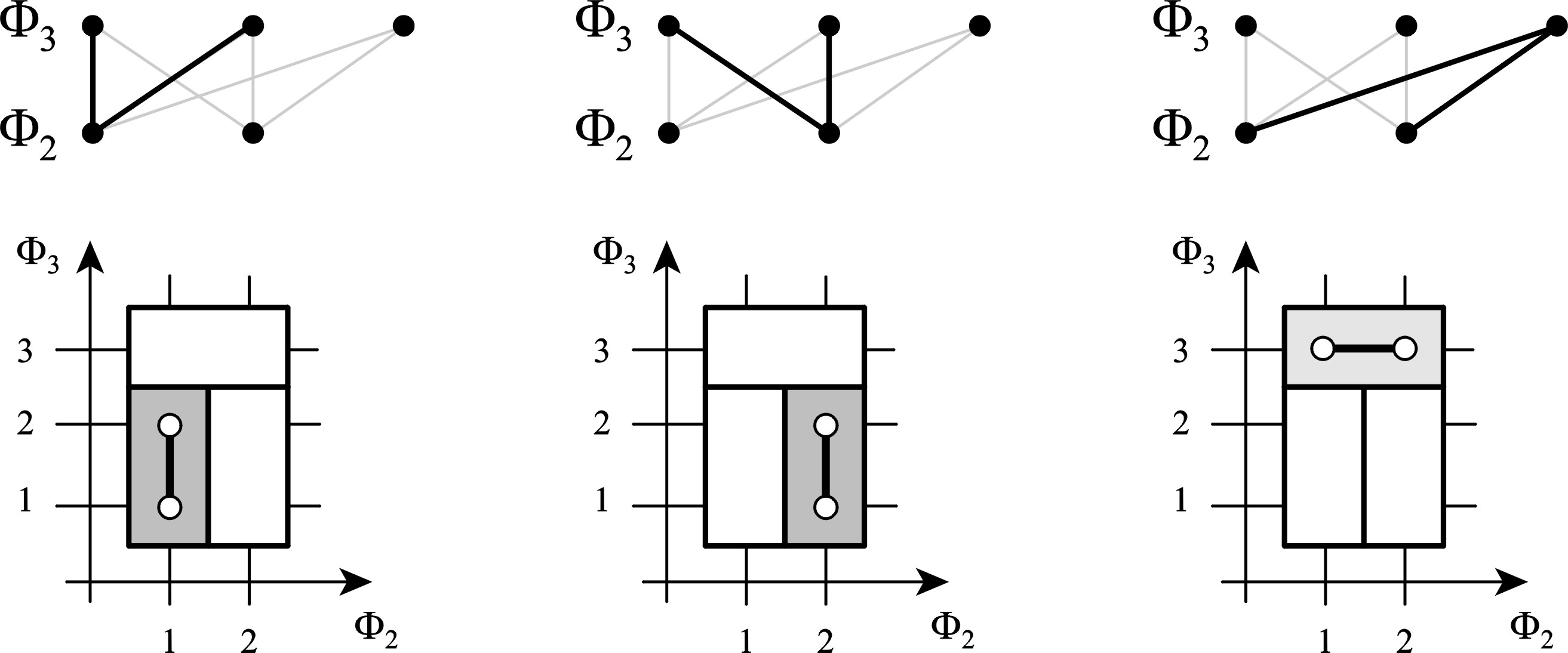}
	\caption{Correspondence between tiling of $F$-box $V_{3,4}$ with two white-dots boxes and tiling of the $\layer{3}{4}$ with two-chain subposets.}
\end{center}
\end{figure}
\noindent \textbf{Recall the fact} (\cite[2007]{11}, \cite[2009]{8}): \textit{Let} $F$\textit{ be an admissible sequence}. \textit{Take any natural numbers} $n,m$ such that $n\geq m$, \textit{then the value of} $F$-\textit{binomial coefficient} $\fnomial{n}{k}$ i\textit{s equal to the number} of sub-boxes \textit{that constitute a}  $\kappa$-\textit{partition of} $m$-\textit{\textit{dimensional}} $F$-\textit{box} $V_{m,n}$ \textit{where} $\kappa = |V_m|$.

\begin{defn}
Let $V_{m,n}$ be a $m$-dimensional $F$-box. Then any $\kappa$-partition into sub-boxes of the form $V_m$ is called tiling of $V_{m,n}$.
\end{defn}

\noindent Hence \textcolor{blue}{\textbf{only these}} partitions of $m$-dimensional box $V_{m,n}$  are admitted  for which all sub-boxes \textbf{are of the form} $V_m$ i.e. we have a kind of  \textcolor{blue}{\textbf{self-similarity}}.\\

\noindent It was shown in \cite[2008]{13} by Maciej Dziemia\'nczuk  that the only admissibility condition  is not sufficient for the existence a tiling for any given $m$-dimensional box $V_{k,n}$. Kwa\'sniewski in \cite[2007]{11} and \cite[2009]{8}  posed the question called \emph{Cobweb Tiling Problem} which we repeat here.

\vspace{0.1cm}

\noindent \textbf{Tiling problem}\\
\noindent Suppose that $F$ is an admissible sequence. Under which conditions any $F$-box $V_{m,n}$ designated by sequence $F$ has a tiling? Find effective characterizations and/or find an algorithm to produce these tilings.

\vspace{0.1cm}

\noindent In \cite[2009]{48} by Maciej Dziemia\'nczuk one  proves the existence of such tiling for certain sub-family of admissible sequences $F$. These  include among others $F=$ Natural numbers, Fibonacci numbers, or $F =  \left\langle n_q \right\rangle_{n\geq 0}$ Gaussian sequence. Original extension of the above  tiling problem onto the general case multi $F$-multinomial coefficients is  proposed  in \cite[2009]{48} , too. Moreover - a reformulation of the present cobweb tiling problem into a clique problem of a graph specially invented for that purpose - is invented.


\vspace{0.3cm}

\noindent Listing. \textbf{15.}  \textit{Fences and \textcolor{blue}{zigzaged roots}}.

\vspace{0.1cm}
\noindent The tiling interpretation of Fibonacci numbers is  known since decades; note for example picture-symbolic language in Chapter 7 from \cite[1989]{GKP 1989}  by Ronald Graham, Donald Ervin  Knuth  and Oren Patashnik. See this picturesque symbolic technique in relevant creative application by Dominique Foata  and  Guo Niu Han in \cite[1994]{Foata 1994 Pisa}.

\vspace{0.1cm}

\noindent  For  \textit{example}  \textbf{0}:  
$F_{n+1}$   equals to the number of compositions of $n$  into \textcolor{red}{\textbf{1}}'s   and  \textcolor{blue}{\textbf{2}}'s; see  \cite[1969]{HoggLind 1969}  by Verner Emil Hoggatt, Jr. and  D.A. Lind  and for more \cite[1975]{AlladiHogg 1975} by Krithnaswami Alladi and  Verner Emil  Hoggatt, Jr.  while the recent generalization is to be found in \cite[2010]{Janjic 2010} by Milan Janji\'c, where the number of all generalized compositions of a natural number is a weighted $r$-generalized Fibonacci number introduced in \cite[2008]{Schork 2008} by Matthias Schork.

\vspace{0.1cm}

\noindent For \textit{example} \textbf{1}:  
$F_{n+1}$   equals to the number of ways  to cover a  $2 \times n$  checkerboard  with  $2 \times 1$  dominoes. There exist also - apart from Listing. \textbf{14.} - at least a twenty years young poset interpretation.

\vspace{0.1cm}

\noindent For \textit{example} \textbf{2}: the $n-th$ Fibonacci number, $F_n$ may be interpreted as the number of ideals in a fence poset as exploited in \cite[1990]{Beck 1990} by Istv\'an Beck. For further continuation of relevant investigation of fences and crowns see: \cite[1991]{ad Beck 1991},  \cite[1992]{ad Beck 1992}, \cite[1995]{ad Beck 1995}  by  Jonathan David Farley,  \cite[2002]{ad Beck 2002}, \cite[2003]{ad Beck 2003}, \cite[2004]{ad Beck 2004}  by Emanuelle Munarini  and  Norma Zagaglia Salvi, \cite[2008]{ad Beck 2008}  by Rodolfo Salvi and  Norma Zagaglia Salvi,
\cite[2009]{Conflitti 2009}  by Alessandro Conflitti.\\

\noindent  - What - if extended? The tiling interpretation of Fibonacci-like numbers is known since several years - up to the knowledge of the present author. 

\vspace{0.2cm}

\noindent For \textit{example} \textbf{3}: according to  the tilings' Combinatorial Theorem 5, p.36  in \cite[2003]{BQ2003}   by Arthur T. Benjamin and  Jennifer J. Quinn we 	ascertain that  for  $H_n = U_n$  the number $s$ from this note recurrence (2)  is interpreted as  equal to the number of colors of squares and  $t$ from this very recurrence (2)  equals to the number of colors of  dominos  while  $H_n = U_{n+1}$    counts colored tilings of length $n$  with squares and dominos. \\ 

\noindent The method of proving various identities in tiling language with the\\ 
$\left\{\textcolor{red}{\textbf{mononimo}},\textcolor{blue}{\textbf{domino}},\right\}$-alphabet ($\equiv$ Morse words alphabet) for number sequences was developed also by Arthur T. Benjamin and  Jennifer J. Quinn and convoked partly in \cite[2003]{BQ2003}. It is becoming now a quite popular tool; for instance see  \cite[2008]{Shattuck 2008}  by Mark Shattuck  or  \cite[2009]{Katz 2009}  by Matt Katz and Catherine Stenson. 

\vspace{0.1cm}

\noindent It is worthy  to ascertain here the pleiad  of investigations and enjoyment with \textit{\textbf{polyominoes}} (P\'olya-ominoes ?) and  P\'olya Festoons - so named in \cite[1991]{Flajolet 1991} by Philippe  Flajolet all these to be placed side by side  with the source idea from \cite[1969]{P69}  by  Gy{\"{o}}rgy  P\'olya  (see more: \textbf{[Ex. q* ; 2]} and \textbf{[Ex. q* ; 3]} below  in subsection \textbf{4.3.}).

\vspace{0.1cm}
\noindent Just to make  a glimpse at this  pleiad  we evoke some references - being totally incomplete  with this undertaking.  Here these are: \cite[1981]{Chvatal 1981}  by C. Berge, C. C. Chen, Vasek Chv\'atal  and S. C. Seow ,  \cite [1990]{Conway Lacarias 1990} by J. H. Conway and  J. C. Lacarias, \cite[1991]{Maylis Delest 1991} by Maylis Delest,  \cite[1991]{Maylis Delest2 1991} by Maylis Delest,  \cite[1995]{Maylis Delest3 1995} by Maylis Delest,  J. P. Dubernard  and I. Dutour,  \cite[1999]{Hickerson 1999}  by   Dean Hickerson ,  \cite[1999]{Renzo Pinzani 1999}  by E. Barcucci, A. Del Lungo, E. Pergola and Renzo Pinzani, \cite[2003]{Del Lungo et all 2003}  by  A. Del Lungo, E. Duchi, A. Frosini, S. Rinaldi;    see also references in  $http://mathworld.wolfram.com/PolyominoTiling.html$  and  $http://mathworld.wolfram.com/Polyomino.html$

\noindent Let us also note  numerous papers by Svjetlan Fereti\'c  et all.  just to mention one \cite[2004]{Feretic 2004}  by  Svjetlan Fereti\'c - directly referring to  Gy{\"{o}}rgy  P\'olya  [quote] \textit{In 1938, P\'olya stated an identity involving the perimeter and area generating function for parallelogram polyominoes. To obtain that identity, P\'olya presumably considered festoons. A festoon (so named by Flajolet) is a closed path w which can be written as w = uv, where each step of u is either $(1, 0)$ or $(0, 1)$, and each step of v is either $(-1, 0)$ or $(0, -1)$.  }
[end of quote] -  ( see \textbf{[Ex. q* ; 2]} and \textbf{[Ex. q* ; 3]} below  in subsection \textbf{4.3.}).

\begin{figure}[ht]
\begin{center}
	\includegraphics[width=85mm]{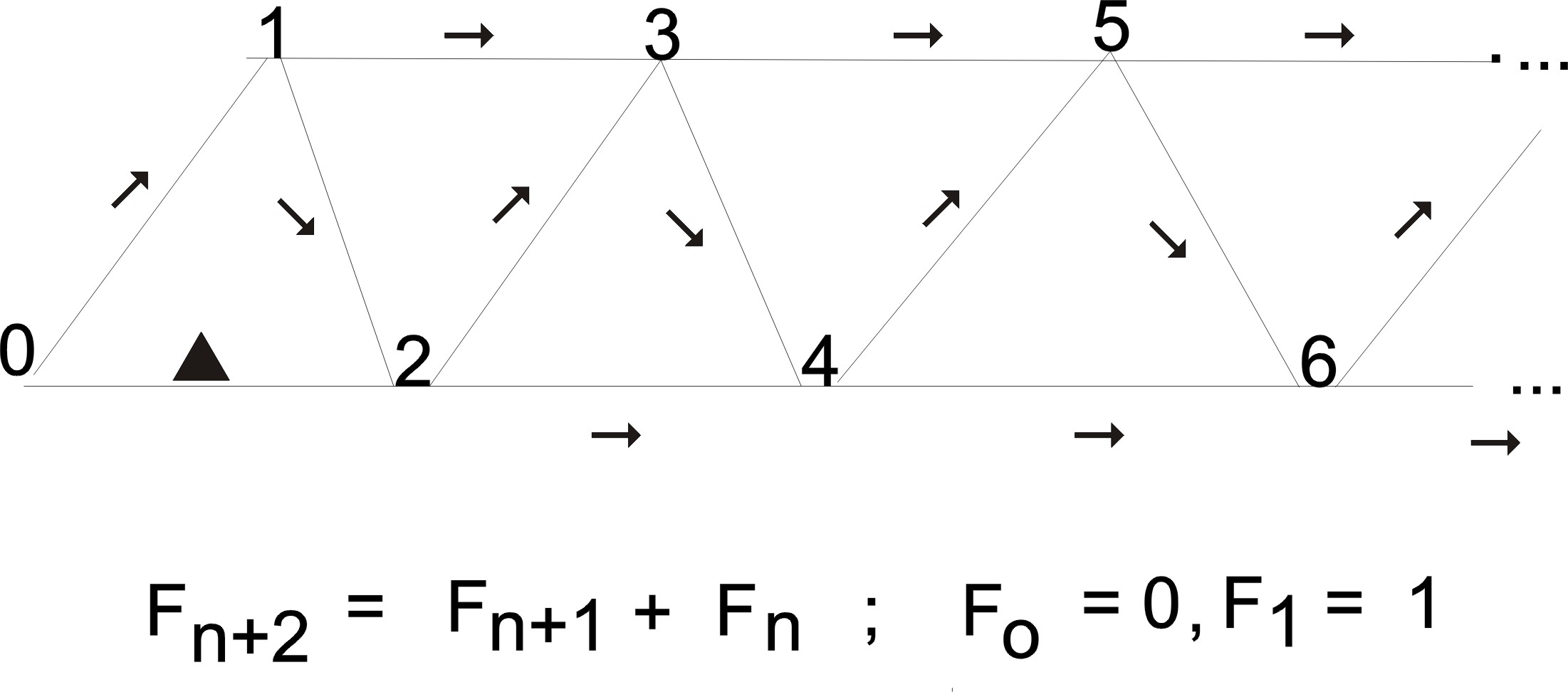}
	\caption {Here are - back to Home \textbf{n} from the Tavern  \textbf{0} - possible night paths - marked with shining arrows.}
\end{center}
\end{figure}

\noindent For \textit{example} \textcolor{blue}{\textbf{4}}: the $n-th$ Fibonacci number $F_n$ may be interpreted as the number of ways to get back Home \textbf{n} from the Tavern  \textbf{0} along paths - \textit{marked with arrows} as displayed by the self-explanatory Fig.6.  (Attention: by paths we mean here paths - \textit{marked with arrows} i.e at least one arrow.)

\vspace{0.1cm}

\noindent For \textit{example} \textcolor{blue}{\textbf{5}}: the $n-th$ \textit{weighted} Fibonacci number i.e.  $F_n\equiv U_n(s,t)$ may be interpreted as the number of  - back to Home \textbf{n} from the Tavern  \textbf{0} - possible weighted  paths - marked with red \textcolor{red}{\textbf{t}} weighted horizontal arrows and with blue  \textcolor{blue}{\textbf{s}} weighted sloping arrows  as displayed by the self-explanatory Fig.7.  \textbf{Attention}:  by paths we mean here paths - \textit{marked with arrows} i.e at least with one arrow In this way there is  no way out, there is no path in case  Home $=$ Tavern i.e. there is no problem in reasonable case of Home $\neq $ Tavern.

\begin{figure}[ht]
\begin{center}
	\includegraphics[width=80mm]{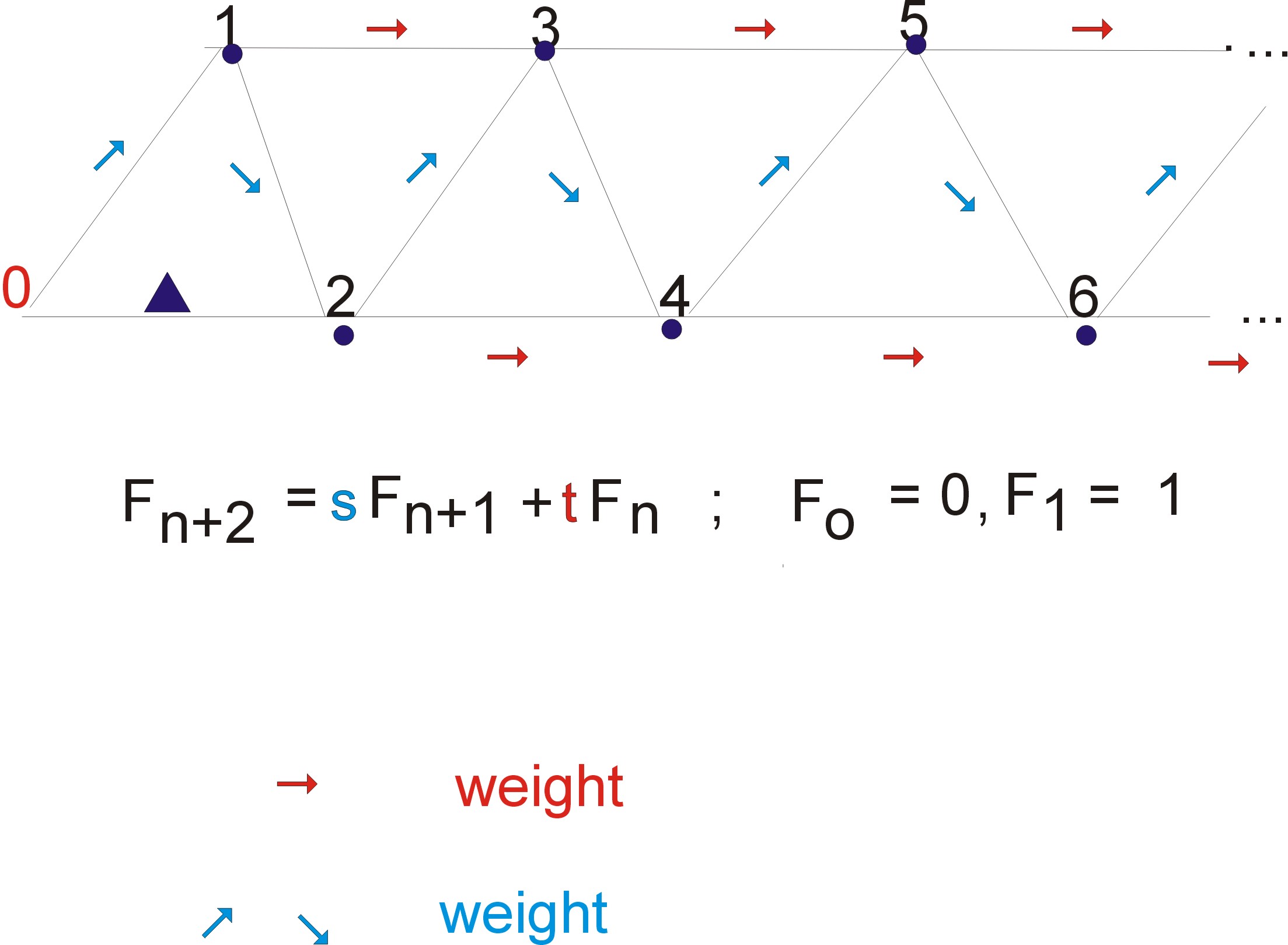}
	\caption {Here are - back to Home \textbf{n} from the Tavern  \textbf{0} - possible weighted  paths - marked with red \textcolor{red}{\textbf{t}} weighted horizontal arrows and with blue  
\textcolor{blue}{\textbf{s}} weighted sloping arrows.}
\end{center}
\end{figure}

\begin{figure}[ht]
\begin{center}
	\includegraphics[width=75mm]{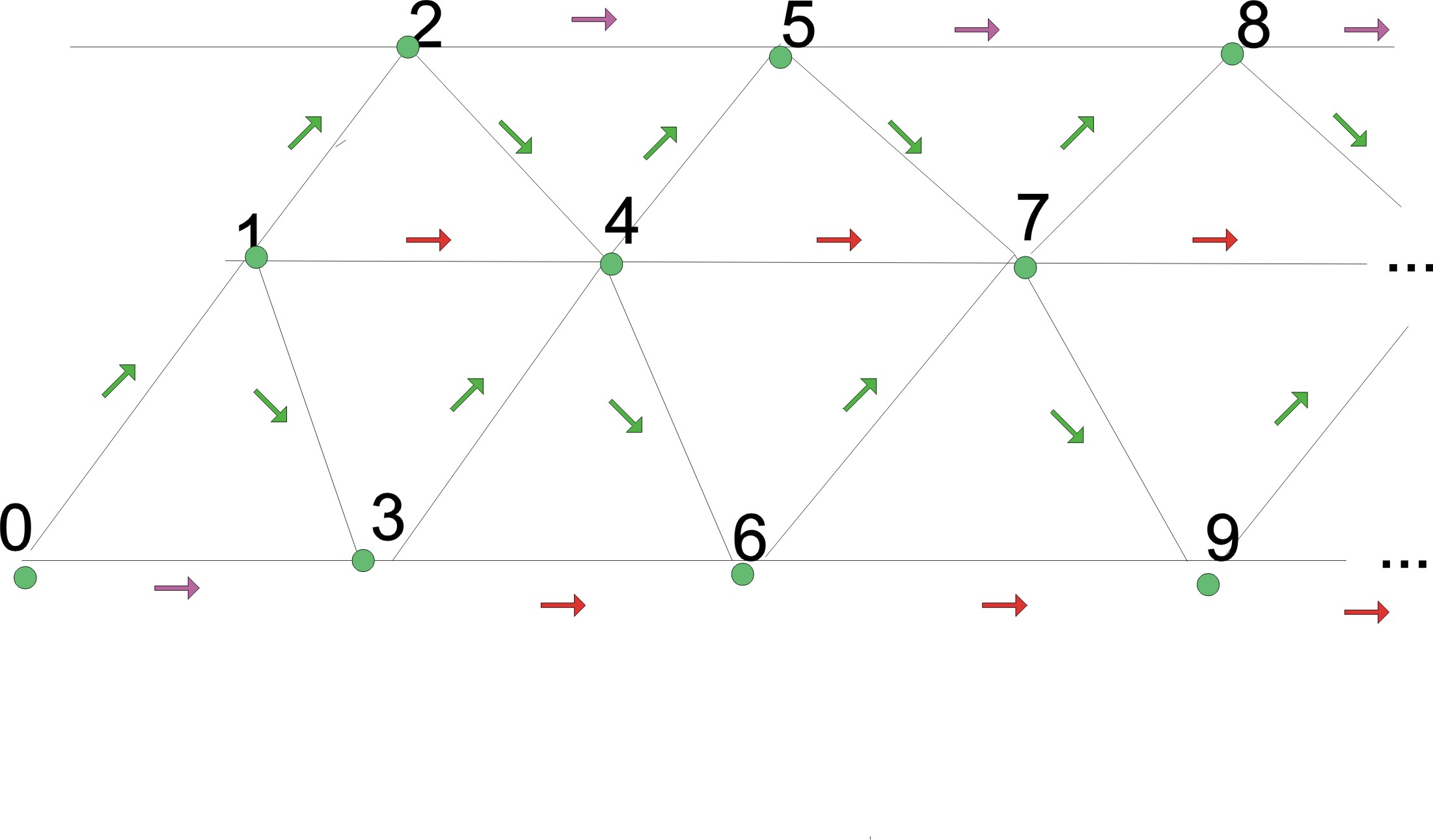}
	\caption {Here are - back to Home \textbf{n} from the Tavern  \textbf{0} - possible two strip highway paths - marked with guiding horizontal and  sloping arrows.}
\end{center}
\end{figure}

\vspace{0.1cm}

\noindent For \textit{example} \textcolor{blue}{\textbf{6}}: the $n-th$ \textit{two strip highway} Fibonacci-like number from the self-explanatory  Fig.8  may be interpreted as the number of  - back to Home \textbf{n} from the Tavern  \textbf{0} - possible  paths - marked with horizontal and  sloping arrows  as displayed by the self-explanatory Fig.8.  

\vspace{0.1cm}

\noindent  According to the Fig.8  we have in this case the following recurrent relations:  $F_0=0$,  $F_1=1$,  $F_2=2$ and  $F_{n+3}= F_{n+2} + F_{n+1}+ F_n$ for $n =1+ k\cdot 3,\,   k\geq 0$, 
$F_{n+3}= F_{n+1}+ F_n$ for $n =0 + k\cdot 3,\:   k\geq 0$  and  $F_{n+3}= F_{n+2} +  F_n$ for $n =2+ k\cdot 3,\:   k\geq 0$. 

\vspace{0.1cm}

\noindent Equivalently: $F_{n+3}= F_{n+2} + F_{n+1}+ F_n = 1 \:mod \ \!3$ for $n =1+ k\cdot 3,\,   k\geq 0$, 
$F_{n+3}= F_{n+1}+ F_n = 0 \:mod \ \!3$ for $n =0 + k\cdot 3,\:   k\geq 0$  and  $F_{n+3}= F_{n+2} +  F_n = 2 \:mod \ \!3 $ for $n =2+ k\cdot 3,\:   k\geq 0$.

\noindent \textbf{Tiling-combinatorial Interpretation of $U$-binomials Question.} 

\vspace{0.3cm}

\noindent What is then combinatorial interpretation of $U$-binomial coefficients stemming from the above different languages examples? Does notation of  items' representing - matter?

\vspace{0.1cm}

\noindent Are the interpretations proposed in \cite[2009]{BSCS} by  Bruce E. Sagan  and Carla  D.  Savage (Listing. \textbf{12.})  and/or  in \cite[2010]{Cigler q-Lucas 2010} by Johann Cigler (Listing. \textbf{11.}) the simplest answer to  \textit{Tiling-combinatorial Interpretation of $U$-binomials Question}?

\vspace{0.2cm}

\noindent  As a \textcolor{blue}{\textbf{Motto}} for an experienced answer we propose:

\vspace{0.1cm}

\noindent \textsl{Science is a language}  (from \cite[1996]{Doron ZeilBerger 1996}   by  Doron Zeilberger); 
     
\vspace{0.1cm}

\noindent \textit{Mathematical notation evolves like all languages}  (from \cite[1992]{Knuth 1992} by Donald Ervin Knuth).

\vspace{0.1cm}

\noindent Another attempt related to  \textit{example} \textbf{3} above by  by  Arthur T. Benjamin and  Sean S. Plott in \cite[2010]{BP}  is presented in the subsequent subsection \textbf{4.2} .


\vspace{0.5cm}

\noindent \textbf{The above list is open and far from complete.}


\vspace{3.cm}

\noindent \textcolor{blue}{\textbf{4.2.} }

\noindent Nevertheless, to this end \textbf{let us discern}  in part - via indicative information - a part of \textbf{Arthur T. Benjamin and coworkers}' recent contribution to the tiling interpretation domain. 
\textbf{Firstly},  let us track in \cite[2003]{BQ2003} the tilings' Combinatorial Theorem 5, p.36. There  for  $H_n = U_n$  the number $s$ from this note recurrence (2)  is interpreted as  equal to the number of colors of squares and  $t$ from this very recurrence (2)  equals to the number of colors of  dominos  while  $H_n = U_{n+1}$    counts colored tilings of length $n$  with squares and dominos.  
\noindent Similarly - also in  \cite[2003]{BQ2003} by Arthur T. Benjamin and  Jennifer J. Quinn  let us consider  the  tilings'  Combinatorial Theorem 6 , p.36.  Here  for $H_n = V_n$  the number   $s$ from this note recurrence (2)  should equal to the number of colors of a square and  $t$ from this very recurrence (2)  equals to the number of colors of a domino while  $H_n = V_n$  counts colored bracelets of length $n$  tiled with squares and dominos. 
\noindent  Bruce E. Sagan  and Carla  D.  Savage in  \cite[2009]{BSCS} refer to well known recurrences:  Identity 73 on p. 38 in \cite{BQ2003} - for (4) in \cite{BSCS} and Identity 94 p. 46 in \cite{BQ2003}  for (5) in \cite{BSCS}. Both (4) and (5) recurrences in \cite[2009]{BSCS} by Bruce E. Sagan  and Carla  D.  Savage  have been evoked  in the illustrative Example 3. Section 3. above.

\vspace{0.1cm}

\noindent \textbf{Firstly}, the paper \cite[2009]{BP}   by Arthur T. Benjamin and  Sean S. Plott  referring to  \cite[2003]{BQ2003}  by Arthur T. Benjamin and Jennifer J. Quinn  should  be notified as it proposes a new formula for $U$-binomials derived via tiling argumentation. Being occasionally nominated  by Arthur T. Benjamin and  Sean S. Plott in  errata \cite[2010]{BP} to this paper  the present author feels entitled to remark also on this errata.

\vspace{0.2cm}

\noindent  According to errata \cite[2010]{BP} by  Arthur T. Benjamin and  Sean S. Plott [quote] " \textsl{The formula for}  ${n \choose k}_F$  \textit{should be multiplied by a factor of} $F_{n - x_k}$, \textsl{which accounts for the one remaining tiling that follows the $f_0$ tiling. Likewise, the formula for  ${n \choose k}_F$ should be multiplied by} $U_{ n - x_k}$."  Our remark is that this errata is unsuccessful. If we follow this errata  then ($x_{k-1} < x_k$) we would have: 

\begin{equation}
	\fnomialF{n}{k}{errata} =
	\sum_{1\leq x_1<x_2<\cdots<x_{k-1}\leq n-1}
	\prod_{i=1}^{k-1} 
	F_{k-i}^{x_i - x_{i-1}-1}
	F_{n-x_i - (k-i) + 1}
	F_{n - x_k},
\end{equation}
where   $F_0 = 0$ and $x_0 = 0$. But the formula (48) implies  for example

$$ 15 = {5 \choose 3}_F \neq \fnomialF{5}{3}{errata} = 11.$$

\vspace{0.1cm}

\noindent The task of finding the correct formula - due to the present author became two months ago an errand - exercise for Maciej Dziemia\'nczuk, a doctoral student from Gda\'nsk University in Poland.
The result - to be quoted below as MD formula (50)  - is his discovery, first announced in the form of a feedback private communication  to the present author: (M. Dziemia\'nczuk  on Mon, Oct 18, 2010 at 6:26 PM)  however still not announced in public.  

\vspace{0.1cm}

\noindent The source of an error in errata is that  ${n \choose k}_F$  should be multiplied \textbf{not by} the factor of $F_{n - x_k}$  \textit{but by} the factor  $F_{n - x_k +1}\equiv f_{n - x_k}$. Then we have

$$
	\fnomialF{n}{k}{now} =
	\sum_{1\leq x_1<x_2<\cdots<x_{k-1}\leq n-1}
	\prod_{i=1}^{k-1} 
	F_{k-i}^{x_i - x_{i-1}-1}
	F_{n-x_i - (k-i) + 1}
	F_{n - x_k +1},
$$
Due to $x_{k-1} < x_k$  the above  formula is equivalent to 

\begin{equation}
	\fnomialF{n}{k}{now} =
	\sum_{1\leq x_1<x_2<\cdots<x_{k-1<x_k}\leq n}
	\prod_{i=1}^{k-1} 
	F_{k-i}^{x_i - x_{i-1}-1}
	F_{n-x_i - (k-i) + 1}
	F_{n - x_k},
\end{equation}
\vspace{0.1cm}

\noindent and this in turn is evidently equivalent to the MD-formula (50) below i.e. (49) is equivalent to the corrected by Maciej Dziemia\'nczuk   Benjamin and Plott  formula from The Fibonacci Quarterly $46/47.1$ (2008/2009), 7-9. 

\vspace{0.2cm}

\noindent Finally here now  MD-formula follows:

\begin{equation}
	\fnomial{n}{k} =
	\sum_{1\leq x_1<x_2<\cdots<x_k\leq n}
	\prod_{i=1}^k 
	F_{k-i}^{x_i - x_{i-1}-1}
	F_{n-x_i - (k-i) + 1},
\end{equation}
where $F_0 = 0$ and $x_0 = 0$. 

\vspace{0.2cm}

\noindent Collaterally  Maciej Dziemia\'nczuk  supplies correspondingly correct formula for  Lucas $U$ - binomial  coefficients $\fnomialF{n}{k}{U}$ :

\begin{align}
	\fnomialF{n}{k}{U} &=
	\sum_{{1\leq x_1<x_2<\cdots<x_{k-1}\leq n} \atop {x_k = x_{k-1}+1}}
	\!\!\!\!\!\!\!\!\!\!\!\!
	s^{x_k - k}
	\left(\prod_{i=1}^{k-1} 
	U_{k-i}^{x_i - x_{i-1}-1}
	U_{n-x_i - (k-i) + 1}\right)
U_{n - x_{k} + 1}
	\\
	&= 
	\sum_{1\leq x_1<x_2<\cdots<x_k\leq n}
	\!\!\!\!\!\!\!\!
s^{x_k - k}
	\prod_{i=1}^k 
	U_{k-i}^{x_i - x_{i-1}-1}
	U_{n-x_i - (k-i) + 1},
\end{align}
where $U_0^t = 0^t = \delta_{t,0}$.


\vspace{0.7cm}

\noindent  \textcolor{red}{\textbf{4.3.}} 
\noindent \textbf{$p,q$-binomials versus \textcolor{blue}{$q*$-binomials combinatorial interpretation},} where   $q* = \frac{p}{q}$ if  $q\neq 0$.

\vspace{0.2cm}

\noindent In the first instance let us  once for all switch off the uninspired  $p\cdot q = 0$ case. Then obligatorily either $q\neq 0$ or  $q\neq 0$. Let then  $q* = \frac{p}{q}$. In this nontrivial case

\begin{equation}
\fnomialF{n}{k}{p,q} =  q^{k(n-k)}\cdot\fnomialF{n}{k}{q*}.
\end{equation}
\vspace{0.2cm}

\noindent Referring to the factor $ q^{k(n-k)}$ as a kind of \textbf{weight},  one may transfer combinatorial interpretation statements on $q*$ binomials $\fnomialF{n}{k}{q*}$ onto 
combinatorial interpretation statements on $p,q$ binomials $\fnomialF{n}{k}{p,q}$ through the agency of (49). Thence , apart from specific  combinatorial interpretations uncovered for the class or subclasses of $p,q$-binomials there might be admitted and respected the "$q*$-overall" combinatorial interpretations transfered from $1,q*$-binomials i.e. from $q*$-binomials onto $p,q$-binomials.

\vspace{0.2cm}

\noindent By no means pretending to be the complete list here comes the skeletonized list of \textbf{[Ex. q* ; k]} examples,  $k\geq 1$.

\vspace{0.4cm}

\noindent \textbf{[Ex. q* ; 1]}  
\noindent The $q*$-binomial coefficient ${m+n \choose m,n}_{q*}$  may be interpreted as a polynomial in $q*$ whose $q*^k$-th coefficient   counts the number of distinct partitions of $k$ elements which fit inside an $m \times n$ rectangle - see  \cite[1976]{Andrews}  by   George Eyre Andrews.

\vspace{0.3cm}

\noindent \textbf{On lattice paths' techniques  - Historical Remark.} It seams to  be desirable  now to  quote here information from \cite[2010]{Katherine 2010}  by Katherine Humphreys based on
\cite[1878]{Whitworth 1878}  by William Allen Whitworth:

\begin{quot}
We find lattice path techniques as early as 1878 in Whitworth to help picture a combinatorial problem, but it is not until the early 1960's that we find lattice path enumeration presented as a mathematical topic on its own. The number of papers pertaining to lattice path enumeration has more than doubled each decade since 1960. 
\end{quot}

\vspace{0.3cm}

\noindent \textbf{[Ex. q* ; 2]}\\ 
\noindent The [Ex. q* ; 2] may be now compiled  with [Ex. q* ; 1] above. For that to do recall that zigzag path is the shortest path  that  starts at
$A = (0,0)$  and ends in $B =  (k,n-k)$ of the $k \times (n-k)$ rectangle; see:  \cite[1962]{P62}  by  Gy{\"{o}}rgy P\'olya   [pp.  68-75],   \cite[1969]{P62}  by  Gy{\"{o}}rgy P\'olya  and
\cite{PA 1971} by  Gy{\"{o}}rgy P\'olya    and  G. L.  Alexanderson.

\vspace{0.1cm}

\noindent Let then     $A_{n,k,\alpha}= $ the number of those $(0,0) \longrightarrow  (k,n-k)$ zigzag paths the area under which  is $\alpha$.

\vspace{0.1cm}

\noindent In  \cite[1969]{P69}  Gy{\"{o}}rgy  P\'olya  using recursion for  $q*$-binomial coefficients proved  that  

$$  \fnomialF{n}{k}{q*} = \sum_{\alpha = 0}^{k(n-k)} A_{n,k,\alpha}\cdot q*^{\alpha}.$$    
from where Gy{\"{o}}rgy  P\'olya infers the following Lemma (\cite[1969]{P69}, p.105) which is named Theorem (p. 104) in more detailed paper \cite[1971]{PA 1971} by  Gy{\"{o}}rgy  P\'olya  and  G. L.  Alexanderson.

\begin{quot}
The number of those zigzag paths the area under which  is $\alpha$  equals  $A_{n,k,\alpha}$.
\end{quot}

\vspace{0.3cm}

\noindent \textbf{[Ex. q* ; 3]}\\ 
\noindent The [Ex. q* ; 3] may be now compared  with [Ex. q* ; 1]. The combinatorial interpretation of ${r+s \choose r,s}_{q*}$ from  [Ex. q* ; 1] had been derived (pp. 106-107) in \cite[1971]{PA 1971} by  Gy{\"{o}}rgy  P\'olya  and  G. L.  Alexanderson, from where - with advocacy from \cite[1971]{Knuth1971}   by  Donald  Ervin  Knuth - we quote the result.

\vspace{0.4cm}

\noindent \textbf{(1971)}: \; ${r+s\choose r,s}_{q*}=$ \textit{ordinary generating function in $\alpha$  powers  of  $q*$  for partitions of $\alpha$ into exactly $r$  non-negative integers none of which exceeds $s$} ,  
\vspace{0.2cm}

\noindent as derived  in \cite[1971]{PA 1971} by  Gy{\"{o}}rgy  P\'olya  and  G. L.  Alexanderson  - see formula (6.9) in \cite{PA 1971}.

\vspace{0.4cm}

\noindent \textbf{(1882)}:\; $\fnomialF{n}{k}{q*} = $ \textit{ordinary generating function in $\alpha$  powers  of  $q$  for partitions of $\alpha$ into at most $k$ parts not exceeding $(n-k)$} ,  

\vspace{0.2cm}

\noindent as recalled in \cite[1971]{Knuth1971} by  Donald  Ervin  Knuth and proved combinatorially in \cite[1882]{SJJ 1882} by James Joseph Sylvester.

\vspace{0.2cm}

\noindent  Let nonce : \: $r+s=n$, $r=k$ then $\textbf{(1971)}\equiv \textbf{(1882)}$ are equal due to    

\begin{equation}
\fnomialF{n}{k}{q*} = \sum_{\alpha = 0}^{k(n-k)} A_{n,k,\alpha}\cdot q*^{\alpha} = \sum_{\alpha = 0}^{r\cdot s} A_{r+s,r,\alpha}\cdot q*^{\alpha} = {r+s \choose r,s}_{q*}.
\end{equation}
where for commodity of comparison formulas in two notations from two papers -  we have been using contractually for a while: $r+s=n$, $r=k$ identifications.

\vspace{0.3cm}

\noindent \textbf{[Ex. q* ; 4]}\\  
\noindent The following was proved in \cite[1961]{Kendall 1961} by Maurice  George Kendall and  Alan Stuart (see p.479 and p.964) and n \cite[1971]{PA 1971} by  Gy{\"{o}}rgy  P\'olya  and  Gerald L.  Alexanderson  (p.106).

\vspace{0.3cm}

\noindent \textit{The area under the zigzag path $=$  The number of inversions  in  the very zigzag path coding sequence.}

\vspace{0.2cm}

\noindent The possible extension of the above combinatorial interpretation onto three dimensional zigzag paths via  "\textit{three-nomials}" was briefly mentioned in \cite[1971]{PA 1971} - see p.108.

\vspace{0.3cm}

\noindent \textbf{[Ex. q* ; 5]}\\  
\noindent The well known (in consequence - \textit{finite geometries'}) interpretation of $\fnomialF{n}{k}{q*}$ coefficient due to Jay Goldman and  Gian-Carlo Rota from \cite[1970]{Rota Goldman 1970} is now  worthy of being recalled;  see also  \cite[1971]{Knuth1971} by  Donald  Ervin  Knuth.

\vspace{0.1cm}

\noindent Let   $V_n$ be an $n$-dimensional vector space over a finite field of $q*$ elements. Then 

\vspace{0.2cm}

\noindent $\fnomialF{n}{k}{q*} $ = \textit{the number of $k$-dimensional subspaces of  $V_n$ }.

\vspace{0.3cm}

\noindent \textbf{[Ex. q* ; 6]}\\  
\noindent This example $=$ the short substantial note  \cite[1971]{Knuth1971} by  Donald  Ervin  Knuth. 

\noindent Compile this example with the example [Ex. q* ; 5] above.

\noindent The essence of a coding  of combinatorial interpretations via bijection between lattices  is the construction of this \textit{coding bijection}  in \cite{Knuth1971}. Namely, let $GF(q*)$ be the Galois field of order $q*$ and let $ V_n \equiv V = GF(q*)^n$ be the $n$-dimensional vector space over $GF(q*)$. Let $[n]= \left\{1,2,...,n\right\}$. Let $\ell(V)$ be the lattice of all subspaces of $ V = GF(q*)^n$  while  $\ell([n])\equiv 2^{[n]}$ denotes the lattice of all subsets of $ [n]$.

\noindent In \cite{Knuth1971} Donald  Ervin  Knuth constructs this  \textit{natural order and rank preserving} map $\Phi$  from the lattice $\ell (V)$ of subspaces  onto the lattice $\ell ([n])\equiv 2^{[n]}$  of subsets of $[n]$. 
  
$$  \ell(V)\;   \stackrel{\Phi}{\rightarrow}\; \ell([n]).$$
\noindent We bethink with some reason whether  this $\Phi$ bijection coding  might be an answer to the subset-subspace problem from subset-subspace problem  from \cite[1998]{Konva 1998} by John Konvalina \textbf{?}

\begin{quot}
 ...the subset-subspace problem (see 6 , 9 , and 3) . The traditional approach to the subset-subspace problem has been to draw the following analogy:   the binomial $\fnomial{n}{k}$  coefficient  counts k-subsets of an n-set, while the analogous Gaussian $\fnomialF{n}{k}{q}$  coefficient counts the number of $k$-dimensional subspaces of an $n$-dimensional finite vector space over the field of 
$q$ elements. 
\noindent The implication from this analogy is that the Gaussian coefficients and related identities tend to the analogous identities for the ordinary binomial coefficients as $q$ approaches 1.
The proofs are often algebraic or mimic subset proofs. But what is the combinatorial reason for the striking parallels between the Gaussian coefficients and the binomial coefficients?
\end{quot}

\noindent According to  Joshef  P. S. Kung  \cite[1995]{Kung 1995} the Knuth's note is not the explanation: 

\begin{quot}
... observation of Knuth yields an order preserving map from $L(V_n(q)$ to Boolean algebra of subsets, but it does not yield a solution to the still unresolved problem of finding a combinatorial interpretation of taking the limit $q \longrightarrow 1$.
\end{quot}Well, perhaps  this limit being performed by $q$-deformed Quantum Mechanics  physicists might be of some help?  There the so called $q$-quantum plain of $q$-commuting variables $x\cdot y- q\cdot y\cdot x =0$   becomes a plane $\mathbb{F}\times \mathbb{F}$  ($\mathbb{F}\  =  \mathbb{R}, \mathbb{C} $,...  $p$-adic fields included) of two commuting variables in the limit $q \longrightarrow 1$. For   see \cite[1953]{ Schutzenberger 1953}  by Marcel-Paul   Sch{\"{u}}tzenberger. For quantum plains - see also  \cite[1995]{Kassel 1995}  by  Christian Kassel. It may deserve notifying that $q$ - extension of of the "classical plane" of commuting variables ($q=1$) seems in a sense ultimate as discussed in \cite[2001]{akk 2001} by  A.K. Kwa\'sniewski

\vspace{0.3cm}

\noindent \textbf{[Ex. q* ; 7]}\\  
\noindent Let us continue the above by further quotation from \cite[1998]{Konva 1998} on generalized binomial coefficients and the subset-subspace problem.

\begin{quot}  
We will show that interpreting the Gaussian coefficients as generalized binomial coefficients of the second kind combinations with repetition reveals the combinatorial connections between not only the binomial coefficients and the Gaussian coefficients, but the Stirling numbers as well. Thus, the ordinary Gaussian coefficient  tends to be an algebraic generalization of the binomial coefficient of the first kind, and a combinatorial generalization of the binomial coefficient of the
second kind.
\end{quot}
Now in order to get more oriented go back to the begining  of  subsection \textbf{ 4.1.} and consult  : Listing. \textbf{1.}, Listing. \textbf{2.}, Listing. \textbf{3.}  which are earlier works and
end up with  \cite[2000]{Konva 2000}  by   John Konvalina on an unified  simultaneous  interpretation of binomial coefficients of both kinds, Stirling numbers of both kinds and Gaussian binomial coefficients  of both kinds.  \textbf{Compare it} then thereafter  with  Listing. \textbf{8.}.

\vspace{0.3cm}




\end{document}